%% file: main.tex
\documentclass[preprint,1p,times,nopreprintline]{elsarticle}

\usepackage{amssymb}
\usepackage{amsmath}
\usepackage{amsthm}

\usepackage{algorithm}
\usepackage{algorithmic}
\usepackage{booktabs}
\usepackage[table]{xcolor}

\usepackage{todonotes}

\newtheorem{theorem}{Theorem}
\newtheorem{definition}[theorem]{Definition}
\newtheorem{example}[theorem]{Example}
\newtheorem{proposition}[theorem]{Proposition}
\newtheorem{remark}[theorem]{Remark}

\journal{Mathematics and Computers in Simulation}

\begin{document}

\begin{frontmatter}



\title{Coordinate-wise splitting algorithms for ODE simulation via Koopman-Lie product formulas} 


\author[lsu]{Arun Banjara} 
\author[lsu]{Ibrahem AlJabea}
\author[ntua]{Theodore Papamarkou}
\author[lsu]{Frank Neubrander}

\affiliation[lsu]{organization={Department of Mathematics, Louisiana State University}, country={USA}}
\affiliation[ntua]{organization={School of Applied Mathematical and Physical Sciences, National Technical University of Athens}, country={Greece}}

\input{abstract.tex}



\begin{keyword}
Operator splitting \sep Koopman semigroup \sep Lie-Trotter \sep Strang splitting \sep product formulas \sep coordinate-wise integration


\end{keyword}

\end{frontmatter}



\input{introduction.tex}
\input{semigroups.tex}
\input{flows.tex}
\input{splitting.tex}
\input{experiments.tex}
\input{conclusion.tex}

\section*{Acknowledgments}

We would like to thank the Department of Mathematics at LSU and the Gordon A.~Cain Center for STEM Literacy for their support of this work.

\appendix

\input{supplementary_algorithms.tex}


\bibliographystyle{elsarticle-num-names} 

\input{references.bbl}

\end{document}

%% file: abstract.tex
\begin{abstract}
We present a computational framework for simulating finite-dimensional ordinary differential equations by combining classical Koopman-Lie product formulas with coordinate-wise frozen subflows. The setting is model-known, since the vector field is assumed to be available, and no data-driven approximation of the Koopman operator is attempted. Under standard assumptions, the Koopman-Lie generator associated with the flow admits a coordinate decomposition into partial generators. This decomposition leads to elementary updates in which all but one state variable are frozen, and the resulting frozen scalar subproblems are evaluated either in closed form or by one-dimensional solves. Lie-Trotter, Strang, and higher-order exponential compositions are then converted into state-update algorithms for two- and three-dimensional systems, with the semigroup and product-formula theory used as background justification for the constructions. We also record the exponential-term counts produced by the recursive constructions used in the implementation. These counts are presented as implementation costs. Numerical experiments on the Lotka-Volterra, Van der Pol, and Lorenz systems compare the coordinate-wise splitting algorithms with high-accuracy RK45 reference solutions using root-mean-square errors and work-precision curves. The results illustrate the practical trade-off between splitting order, number of time steps, number of exponential factors, and runtime.
\end{abstract}

%% file: introduction.tex
\section{Introduction}
\label{sec:introduction}

Let $\Omega\subseteq \mathbb{R}^N$ be a state space and let $F=(F_1,\ldots,F_N):\Omega\to\mathbb{R}^N$ be a vector field. We consider the autonomous initial-value problem
\begin{equation}
\label{IVP}
u'(t)=F(u(t)),
~
u(0)=u_0\in\Omega .
\end{equation}
When \eqref{IVP} generates a flow $\varphi_t:\Omega\to\Omega$, the associated Koopman semigroup acts on observables by $T(t)g(u)=g(\varphi_t(u))$. This gives a linear representation of the nonlinear dynamics at the level of observables. The operator-theoretic viewpoint goes back to~\citet{Koopman2,Koopman1}. Related flow and generator ideas also appear in the work of~\citet{New62,New61}. For jointly continuous flows,~\citet{JD1,JD2} developed a semigroup framework in which such generators can be studied on spaces of bounded continuous functions~\citep{JD1,JD2}.

For smooth finite-dimensional systems, the Koopman generator has the expression
\[
\mathcal{K}g(u)
=
\sum_{i=1}^N F_i(u)\frac{\partial g}{\partial u_i}(u).
\]
This representation suggests a coordinate decomposition of $\mathcal{K}$ into partial generators. Each partial generator corresponds to a frozen-coordinate equation in which all variables except one are held fixed. If the corresponding frozen flow remains in $\Omega$, then this gives an elementary state update. These elementary updates can then be composed by product formulas such as the Lie-Trotter and Strang formulas, or by higher-order exponential compositions~\citep{KHH1,strang1968construction,chorin1978product,JD5,JD6,JD7}.

The present paper deals with this model-known setting. The vector field $F$ is assumed to be available. Alternative approaches construct a finite-dimensional data-driven approximation of the Koopman operator or its generator based on applications of DMD, EDMD, or related dictionary-based methods~\citep{379schmid2010dynamic,366rowley2009spectral,227kutz2016dynamic,EDMD2}. Instead, we use the Koopman-Lie notation to organize a family of coordinate-wise splitting algorithms for state-space simulation. The numerical output is an approximate state map $\varphi_n(t,\cdot)$. Once this state map is computed, any chosen observable can be evaluated by composition with $\varphi_n(t,\cdot)$.

The comparison with Runge-Kutta methods is carried out along the same lines. A high-accuracy RK45 solution provides a reference trajectory for measuring errors in the numerical experiments. The purpose of our splitting scheme is to exploit frozen subflows that may be cheaply solvable. This can be useful when the vector field decomposes into coordinate mechanisms and when higher-order compositions provide a favorable trade-off between time-step count, number of elementary maps, and observed accuracy.

\paragraph{Main contributions.}
The contribution of this paper is computational. We make use of classical Koopman-Lie product formulas in coordinate-wise algorithms for finite-dimensional ordinary differential equation (ODE) simulation. The main implementation step is to replace the full flow by compositions of frozen one-dimensional subflows, evaluated either in closed form or by one-dimensional solves when closed forms are not available. We also record the number of exponential factors produced by the recursive constructions used in our implementation. These counts are implementation costs for the constructions considered here. Finally, we test the resulting Lie-Trotter, Strang, and higher-order schemes on the Lotka-Volterra, Van der Pol, and Lorenz systems using root-mean-square error (RMSE) and work-precision comparisons against high-accuracy RK45 reference solutions.

%% file: semigroups.tex
\section{Strongly continuous and bi-continuous semigroups}
\label{sec:koopmansemigroup}

We recap the semigroup terminology used in this paper. Let $X$ be a complex Banach space with norm $\|\cdot\|$, and let $\mathcal{L}(X)$ denote the bounded linear operators on $X$. A family $\{T(t)\}_{t\ge 0}\subseteq\mathcal{L}(X)$ is a one-parameter semigroup if $T(0)=I,~T(t+s)=T(t)T(s)$ for all $s,t\ge 0$. It is a strongly continuous semigroup if, for every $x\in X$, the map $t\mapsto T(t)x$ is continuous from $[0,\infty)$ into $X$ with respect to the norm topology.

The generator of a strongly continuous semigroup $\{T(t)\}_{t\ge 0}$ is the operator $(\mathcal{A},D(\mathcal{A}))$ defined by
\[
\mathcal{A}x
=
\lim_{t\to 0^+}\frac{T(t)x-x}{t},
\]
on the domain
\[
D(\mathcal{A})
=
\left\{
x\in X:
\lim_{t\to 0^+}\frac{T(t)x-x}{t}
\text{ exists in }X
\right\}.
\]
For strongly continuous semigroups, the generator $\mathcal{A}$ is closed, and its domain $D(\mathcal{A})$ is dense in $X$. Moreover, there exist constants $M\ge 1$ and $\omega\in\mathbb{R}$ such that $\|T(t)\|\le M e^{\omega t},~t\ge 0$. For further details, see, for example,~\citet{nagel2,pazy,gold}.

For Koopman semigroups on spaces of bounded continuous functions, strong continuity in the sup norm can fail. The following example explains why the bi-continuous framework is useful.

\begin{example}
\label{nowhere}
Let $\Omega=[0,\infty)$ and let $\varphi_t(x)=x+t$. This is the flow generated by $x'(t)=1$ with $x(0)=x$. The induced Koopman semigroup is
\[
T(t)g(x)=g(\varphi_t(x))=g(x+t).
\]
On $C_0([0,\infty),\mathbb{C})$, this semigroup is strongly continuous. On $C_b([0,\infty),\mathbb{C})$, it is not strongly continuous in the sup norm. Indeed, for $g(x)=e^{ix^2}$,
\begin{align*}
|T(t+r)g(x)-T(t)g(x)|
&=
|e^{i\phi(x)}-1|,\\
\phi(x)&=(x+t+r)^2-(x+t)^2=2rx+r^2+2tr.
\end{align*}
For $r>0$, the function $\phi$ is affine in $x$ with positive slope, hence attains values congruent to $\pi$ modulo $2\pi$ at suitable $x\ge 0$. At such points $|e^{i\phi(x)}-1|=2$, so $\|T(t+r)g-T(t)g\|_\infty=2$ for every $t\ge 0$ and every $r>0$. Thus, $t\mapsto T(t)g$ is not norm-continuous at any $t\ge 0$.
\end{example}

We now define the locally convex topology used throughout the paper. If $\Omega$ is a locally compact metric space, let $C_b(\Omega)=C_b(\Omega,\mathbb{C})$. We write $\tau=\tau_{\mathrm{co}}$ for the compact-open topology on $C_b(\Omega)$. This topology is generated by the seminorms $p_K(g)=\sup_{x\in K}|g(x)|$, where $K$ ranges over the compact subsets of $\Omega$. Thus convergence in $\tau$ is uniform convergence on compact subsets of $\Omega$. The norm on $C_b(\Omega)$ is always the sup norm $\|g\|_\infty=\sup_{x\in\Omega}|g(x)|$.

The shift semigroup presented in Example~\ref{nowhere} is not strongly continuous on $C_b([0,\infty),\mathbb{C})$ with respect to $\|\cdot\|_\infty$, but it is continuous with respect to $\tau_{\mathrm{co}}$. More specifically, if $g\in C_b([0,\infty),\mathbb{C})$, then $T(t)g\to g$ uniformly on compact subsets as $t\to 0^+$. This motivates the following definition.

\begin{definition}
\label{Def:re-admissable}
Let $(X,\|\cdot\|)$ be a Banach space and let $\tau$ be a locally convex topology on $X$ generated by seminorms $\{p_\alpha\}_{\alpha\in I}$. We call $(X,\|\cdot\|,\tau)$ a bi-admissible Banach space if the following conditions hold:
\begin{itemize}
\item $\|x\|=\sup_{\alpha\in I}p_\alpha(x)$ for all $x\in X$.
\item Every norm-bounded $\tau$-Cauchy sequence in $X$ is $\tau$-convergent in $X$.
\item $x=0$ if and only if $p_\alpha(x)=0$ for all $\alpha\in I$.
\item The dual space $(X,\tau)^*$ is norming for $(X,\|\cdot\|)$, that is,
\[
\|x\|
=
\sup\{|\phi(x)|:\phi\in (X,\tau)^*,\ \|\phi\|\le 1\}.
\]
\end{itemize}
\end{definition}

\begin{definition}
\label{defBi}
Let $(X,\|\cdot\|,\tau)$ be a bi-admissible Banach space. A family $\{T(t)\}_{t\ge 0}\subseteq \mathcal{L}(X)$ is a bi-continuous semigroup of type $\omega$ if the following conditions hold:
\begin{itemize}
\item $T(0)=I$ and $T(t+s)=T(t)T(s)$ for all $s,t\ge 0$.
\item There exist $M\ge 1$ and $\omega\in\mathbb{R}$ such that $\|T(t)\|\le M e^{\omega t},~t\ge 0$.
\item For every $x\in X$, the map $t\mapsto T(t)x$ is $\tau$-continuous on $[0,\infty)$.
\item The family is locally bi-equicontinuous: if $x_n\to 0$ in $\tau$ and $\sup_n\|x_n\|<\infty$, then $T(t)x_n\to 0$ in $\tau$, uniformly for $t$ in compact subintervals of $[0,\infty)$.
\end{itemize}
\end{definition}

Let $\{T(t)\}_{t\ge 0}$ be a bi-continuous semigroup on $(X,\|\cdot\|,\tau)$. Its generator $(\mathcal{K},D(\mathcal{K}))$ is defined by
\[
\mathcal{K}x
=
\tau\text{-}\lim_{t\to 0^+}\frac{T(t)x-x}{t},
\]
with domain
\[
D(\mathcal{K})
=
\left\{
x\in X:
\tau\text{-}\lim_{t\to 0^+}\frac{T(t)x-x}{t}
\text{ exists and }
\sup_{0<t\le 1}
\left\|
\frac{T(t)x-x}{t}
\right\|
<\infty
\right\}.
\]
A subspace $D\subset X$ is called bi-dense if every $x\in X$ is the $\tau$-limit of a norm-bounded sequence from $D$. An operator $(\mathcal{K},D(\mathcal{K}))$ is called bi-closed if, whenever $x_n\in D(\mathcal{K})$, $\sup_n\{\|x_n\|,\|\mathcal{K}x_n\|\}<\infty,~x_n\to x$ in $\tau,~\mathcal{K}x_n\to y$ in $\tau$, then $x\in D(\mathcal{K})$ and $\mathcal{K}x=y$.

For $C_b(\Omega)$ equipped with the compact-open topology $\tau_{\mathrm{co}}$, Koopman semigroups generated by jointly continuous global flows constitute examples of bi-continuous semigroups. This is the setting used in the next section.

%% file: flows.tex
\section{Flows and Koopman semigroups}
\label{sec:flow}

This section recaps on the background needed for the splitting constructions of the present paper. We recall the connection between flows and Koopman semigroups.

Let $\Omega$ be a topological space. Since only nonnegative times are used, we work with semiflows, yet we use the term flow for brevity. A map $\varphi:[0,\infty)\times\Omega\to\Omega,~(t,x)\mapsto \varphi_t(x)$, is an autonomous, $\Omega$-invariant flow if $\varphi_0(x)=x,~\varphi_t(\varphi_s(x))=\varphi_{t+s}(x)$ for all $x\in\Omega$ and all $s,t\ge 0$. The flow is jointly continuous if $\varphi$ is continuous as a map from $[0,\infty)\times\Omega$ into $\Omega$.

We next relate this definition to ODEs. Let $O\subseteq\mathbb{R}^N$ be open, let $\Omega\subseteq O$, and let $F:O\to\mathbb{R}^N$ be a vector field. A curve $u:[0,\infty)\to\Omega$ is a classical solution of $u'(t)=F(u(t)),~u(0)=u_0\in\Omega$, if $u$ is differentiable and satisfies the equation for all $t\ge 0$. If, for every $u_0\in\Omega$, this problem has a unique global solution $u(t;u_0)$ that remains in $\Omega$, and if $(t,u_0)\mapsto u(t;u_0)$ is jointly continuous, then $\varphi_t(u_0)=u(t;u_0)$ defines a jointly continuous, $\Omega$-invariant flow.

Let $C_b(\Omega)=C_b(\Omega,\mathbb{C})$. A flow $\varphi$ induces the Koopman operators $T(t)g(x)=g(\varphi_t(x)),~g\in C_b(\Omega),~x\in\Omega,~t\ge 0$. The operators $T(t)$ are linear. They satisfy $T(0)=I$ and $T(t+s)=T(t)T(s)$. If each $\varphi_t$ is continuous, then $T(t)$ maps $C_b(\Omega)$ into itself. Moreover, $\|T(t)g\|_\infty\le \|g\|_\infty,~t\ge 0$. Thus, the Koopman operators are contractions on $C_b(\Omega)$.

The following proposition gives the flow-to-semigroup implication needed for our construction.

\begin{proposition}
\label{D-N}
Let $\Omega$ be a locally compact metric space, and let $\varphi:[0,\infty)\times\Omega\to\Omega$ be a jointly continuous, $\Omega$-invariant flow. Let $\tau=\tau_{\mathrm{co}}$ be the compact-open topology on $C_b(\Omega)$. Then the Koopman operators $T(t)g=g\circ\varphi_t,~t\ge 0$, form a bi-continuous contraction semigroup on $(C_b(\Omega),\|\cdot\|_\infty,\tau)$.
\end{proposition}

\begin{proof}
The semigroup identities follow from the flow identities. The contraction estimate follows from
\[
\|T(t)g\|_\infty
=
\sup_{x\in\Omega}|g(\varphi_t(x))|
\le
\|g\|_\infty .
\]
It remains to check the continuity properties with respect to $\tau$. Let $K\subset\Omega$ be compact and let $T>0$. Since $\varphi$ is jointly continuous, the set $L=\varphi([0,T]\times K)$ is compact. If $g\in C_b(\Omega)$, then $g$ is uniformly continuous on $L$. Hence, $T(t)g\to g$ uniformly on $K$ as $t\to 0^+$. The same argument, applied to $\varphi_t$ locally in time, gives $\tau$-continuity of $t\mapsto T(t)g$ on $[0,\infty)$.

For local bi-equicontinuity, suppose that $g_n\to 0$ in $\tau$ and that $\sup_n\|g_n\|_\infty<\infty$. Then, for the compact set $L$ defined above,
\[
\sup_{0\le t\le T}\sup_{x\in K}|T(t)g_n(x)|
=
\sup_{0\le t\le T}\sup_{x\in K}|g_n(\varphi_t(x))|
\le
\sup_{y\in L}|g_n(y)|
\to 0.
\]
Thus $T(t)g_n\to 0$ in $\tau$, uniformly for $t\in[0,T]$. This proves bi-continuity.
\end{proof}

The generator of the Koopman semigroup is the bi-continuous generator defined
in Section~\ref{sec:koopmansemigroup}. Thus,
\[
\mathcal{K}g
=
\tau\text{-}\lim_{t\to 0^+}\frac{T(t)g-g}{t},
\]
with domain
\[
D(\mathcal{K})
=
\left\{
g\in C_b(\Omega):
\tau\text{-}\lim_{t\to 0^+}
\frac{T(t)g-g}{t}
\text{ exists, }
\sup_{0<t\le 1}
\left\|
\frac{T(t)g-g}{t}
\right\|_\infty
<\infty
\right\}.
\]
Equivalently, for $g\in D(\mathcal{K})$,
\[
(\mathcal{K}g)(x)
=
\lim_{t\to 0^+}
\frac{g(\varphi_t(x))-g(x)}{t},
~
x\in\Omega.
\]
The full generator domain $D(\mathcal{K})$ in the bi-continuous setting is defined via the pointwise formula, the compact-open limit, and the bounded-difference-quotient condition.

For finite-dimensional ODEs, this generator has the directional-derivative form on smooth observables for which the expression belongs to the generator domain. If $\Omega\subseteq O\subseteq\mathbb{R}^N$, $F=(F_1,\ldots,F_N)$, and the flow is generated by $u'=F(u)$, then
\[
\mathcal{K}g(u)
=
\nabla g(u)\cdot F(u)
=
\sum_{i=1}^N F_i(u)\frac{\partial g}{\partial u_i}(u).
\]
This identity holds, for example, for $g\in C^1(O)$ when $g|_\Omega\in C_b(\Omega)$, the function $u\mapsto \nabla g(u)\cdot F(u)$ belongs to $C_b(\Omega)$, and the limit for $\mathcal{K}$ holds in the compact-open topology with bounded difference quotients. In the algorithms of the present paper, this formula is used as the coordinate representation of the Koopman-Lie generator associated with the known vector field $F$.

\begin{remark}
The Koopman semigroup acts on observables, but the numerical algorithms in this paper are implemented through the associated state maps, without needing to construct a finite-dimensional dictionary approximation of $T(t)$ or $\mathcal{K}$. Once an approximate state map $\varphi_n(t,\cdot)$ has been computed, an observable is evaluated by composition, namely $g(\varphi_n(t,\cdot))$.
\end{remark}

%% file: splitting.tex
\section{Splitting methods}
\label{sec:Splitting methods}

We now describe the coordinate-wise splitting construction used in the numerical algorithms. The construction is stated at the level of state maps and then translated into Koopman notation. The Koopman semigroups act on observables, while the numerical implementation advances points in the state space.

Let $O\subseteq\mathbb{R}^N$ be open, let $\Omega\subseteq O$, and let $F=(F_1,\ldots,F_N):O\to\mathbb{R}^N$ be a vector field. Assume that the initial-value problem $u'(t)=F(u(t)),~u(0)=u_0\in\Omega$, generates a unique global, jointly continuous, $\Omega$-invariant flow $\varphi_t:\Omega\to\Omega$. The corresponding Koopman semigroup is $T(t)g(u)=g(\varphi_t(u)),~g\in C_b(\Omega)$.

For each $i=1,\ldots,N$, define the $i$-th coordinate vector field $F^{[i]}(u)=F_i(u)e_i$, where $e_i$ is the $i$-th coordinate vector in $\mathbb{R}^N$. We assume that each $F^{[i]}$ generates a unique global, jointly continuous, $\Omega$-invariant frozen-coordinate flow $\Phi_i:[0,\infty)\times\Omega\to\Omega$. In coordinates, if $u=(u_1,\ldots,u_N)$, then
\[
\Phi_i(t,u)
=
(u_1,\ldots,u_{i-1},\psi_i(t,u),u_{i+1},\ldots,u_N),
\]
where $\psi_i(t,u)$ solves the scalar frozen problem
\[
w'(t)
=
F_i(u_1,\ldots,u_{i-1},w(t),u_{i+1},\ldots,u_N),
~
w(0)=u_i .
\]
Thus, all coordinates except the $i$-th one are held fixed during the $i$-th frozen substep.

Each frozen flow induces a Koopman semigroup $T_i(t)g(u)=g(\Phi_i(t,u))$, where $g\in C_b(\Omega)$. By Proposition~\ref{D-N}, each $\{T_i(t)\}_{t\ge 0}$ is a bi-continuous contraction semigroup on
$(C_b(\Omega),\|\cdot\|_\infty,\tau_{\mathrm{co}})$. We denote its generator by $\mathcal{K}_i$. The notation $e^{t\mathcal{K}_i}g$ will be used as shorthand for $T_i(t)g$.

The following proposition states the coordinate form of the generators on smooth observables, which is used as the generator decomposition in this paper.

\begin{proposition}
\label{prop:coordinate-generators}
Assume that $F$ generates a unique global, jointly continuous, $\Omega$-invariant flow, and that each frozen vector field $F^{[i]}$ generates a unique global, jointly continuous, $\Omega$-invariant frozen-coordinate flow. Let $g\in C^1(O,\mathbb{C})$, and write $g$ for its restriction to $\Omega$. Suppose that
$
g|_\Omega\in C_b(\Omega),
~
F_i\,\frac{\partial g}{\partial u_i}\Big|_\Omega\in C_b(\Omega)$ for $i=1,\ldots,N$, and $F\cdot\nabla g\big|_\Omega\in C_b(\Omega)$. Then $g$ belongs to the domains of $\mathcal{K}$ and $\mathcal{K}_i$, and
\[
\mathcal{K}_i g(u)
=
F_i(u)\frac{\partial g}{\partial u_i}(u),
~
i=1,\ldots,N.
\]
Moreover,
\[
\mathcal{K}g(u)
=
F(u)\cdot\nabla g(u)
=
\sum_{i=1}^N \mathcal{K}_i g(u).
\]
\end{proposition}

\begin{proof}
We give the argument for $\mathcal{K}_i$; the proof for $\mathcal{K}$ is the same, with $\Phi_i$ replaced by $\varphi$. For $u\in\Omega$, the chain rule gives
\[
\frac{d}{ds}g(\Phi_i(s,u))
=
F_i(\Phi_i(s,u))
\frac{\partial g}{\partial u_i}(\Phi_i(s,u)).
\]
Hence,
\[
\frac{T_i(t)g(u)-g(u)}{t}
=
\frac{1}{t}\int_0^t
F_i(\Phi_i(s,u))
\frac{\partial g}{\partial u_i}(\Phi_i(s,u))\,ds .
\]
Since $\Phi_i$ is jointly continuous, the integrand converges uniformly on compact subsets of $\Omega$ to $F_i(u)\frac{\partial g}{\partial u_i}(u)$ as $t\to 0^+$. The assumed boundedness gives the bounded-difference-quotient condition in the definition of the bi-continuous generator. Therefore,
\[
\mathcal{K}_i g(u)
=
F_i(u)\frac{\partial g}{\partial u_i}(u).
\]
Summing these identities gives the stated formula for $\mathcal{K}g$.
\end{proof}

The product formulas below are applied to the frozen semigroups $T_i(t)$, or equivalently to their associated state maps $\Phi_i(t,\cdot)$. At the state level, the order of composition is read from right to left. For example,
\[
T_1(h)T_2(h)g
=
g\circ\Phi_2(h,\cdot)\circ\Phi_1(h,\cdot).
\]
Thus, the observable product $T_1(h)T_2(h)$ corresponds to first applying the first frozen state map $\Phi_1(h,\cdot)$ and then applying the second frozen state map $\Phi_2(h,\cdot)$.

All theoretical splitting formulas in this paper use exact frozen subflows. In some problems, these maps are available in closed form. In other problems, they may be represented by an implicit quadrature formula or computed by a scalar numerical solution. If the frozen subflows are approximated numerically, the additional local error must be included in the error analysis. The numerical examples in Section~\ref{sec:experiments} use closed-form frozen maps.

\subsection{Exponential product formulas}
\label{subsec:exp_prod_formulas}

The known results in this subsection provide the background related to the product formulas used in this paper. They justify the passage from exact frozen subflows to approximate state maps when the standard assumptions of the relevant product formula are satisfied.

Let $(X,\|\cdot\|,\tau)$ be a bi-admissible Banach space, and let $V:[0,\infty)\to\mathcal{L}(X)$ satisfy $V(0)=I$. A bi-continuous version of Chernoff's product formula states that, under the usual stability and consistency hypotheses, the powers of $V(t/n)$ converge to the semigroup generated by the closure of the derivative of $V$ at zero. We use the following form.

\begin{theorem}[Chernoff product formula in the bi-continuous setting~\citep{JD3}]
\label{CPFBi}
Let $V:[0,\infty)\to\mathcal{L}(X)$ satisfy $V(0)=I$. Assume that there are constants $M\ge 1$ and $\omega\in\mathbb{R}$ such that $\|V(h)^m\|\le M e^{\omega mh}, ~ h\ge 0, ~ m\in\mathbb{N}$. Assume also that, for every $T>0$, the family $\{V(h)^m: h\ge 0,\ m\in\mathbb{N},\ mh\le T\}$ is locally bi-equicontinuous. Let $D\subset X$ be bi-dense, and suppose that, for every $f\in D$, the limit
\[
\mathcal{A}_0 f
=
\tau\text{-}\lim_{h\to 0^+}\frac{V(h)f-f}{h}
\]
exists and satisfies
\[
\sup_{0<h\le 1}
\left\|
\frac{V(h)f-f}{h}
\right\|
<\infty .
\]
If the bi-closure $\mathcal{A}$ of $(\mathcal{A}_0,D)$ generates a bi-continuous semigroup $\{T(t)\}_{t\ge 0}$, then
\[
T(t)f
=
\tau\text{-}\lim_{n\to\infty}
V\!\left(\frac{t}{n}\right)^n f,
\]
for all $f\in X$, uniformly for $t$ in compact subintervals of $[0,\infty)$.
\end{theorem}

In the strongly continuous case, the same statement is read with $\tau$ replaced by the norm topology. In the bi-continuous case, the derivative condition in Theorem~\ref{CPFBi} is a $\tau$-limit together with a bounded-difference-quotient condition. This is the same type of generator condition used in Section~\ref{sec:koopmansemigroup}.

We now apply this notation to Koopman semigroups. Let $\varphi_t$ be the full flow generated by $u'=F(u)$, and let $T(t)g=g\circ\varphi_t$ be the associated Koopman semigroup on $C_b(\Omega)$. Let $\Phi_i(t,\cdot)$ be the frozen state maps introduced above, and write $T_i(t)g=g\circ\Phi_i(t,\cdot)$. A one-step splitting map $\Psi_h:\Omega\to\Omega$ induces an operator $V(h)g=g\circ\Psi_h$. If the hypotheses of Theorem~\ref{CPFBi} hold for this $V$, then
\begin{equation}
\label{Arun}
g(\varphi_t(\cdot))
=
T(t)g
=
\tau\text{-}\lim_{n\to\infty}
V\!\left(\frac{t}{n}\right)^n g .
\end{equation}
Equivalently, if $\Psi_h^{[n]}$ denotes the $n$-fold composition of $\Psi_h$ with itself, then
\[
V\!\left(\frac{t}{n}\right)^n g(u)
=
g\!\left(\Psi_{t/n}^{[n]}(u)\right).
\]
Thus, \eqref{Arun} gives convergence of observables along the approximate state maps.

The following proposition converts this observable convergence into pointwise convergence of state values.

\begin{proposition}
\label{Frank}
Let $\Omega$ be a metric space, and let $u_n,u\in\Omega$. Then $u_n\to u$ in $\Omega$ is equivalent to $g(u_n)\to g(u)$ for every $g\in C_b(\Omega,\mathbb{R})$.
\end{proposition}

\begin{proof}
If $u_n\to u$, then $g(u_n)\to g(u)$ for every continuous $g$. Conversely, suppose that $g(u_n)\to g(u)$ for every $g\in C_b(\Omega,\mathbb{R})$. Let $d$ be the metric on $\Omega$, and define
\[
g_u(v)=\frac{d(v,u)}{1+d(v,u)} .
\]
Then $g_u\in C_b(\Omega,\mathbb{R})$ and $g_u(u)=0$. Hence $g_u(u_n)\to 0$, which implies $d(u_n,u)\to 0$.
\end{proof}

It follows that, for each fixed $u\in\Omega$, convergence in \eqref{Arun} for all bounded continuous observables implies $\Psi_{t/n}^{[n]}(u)\to \varphi_t(u)$. In this way, convergence of the observable product formula implies convergence of the corresponding state-space splitting algorithm.

We next state the two basic real-time splitting formulas. Assume that the relevant product-formula hypotheses hold for the coordinate generators $\mathcal{K}_1,\ldots,\mathcal{K}_N$, and that their sum has bi-closure equal to the Koopman generator $\mathcal{K}$ of the full flow. The Lie--Trotter one-step operator is
\begin{equation}
\label{LTstep}
V_{\mathrm{LT}}(h)
=
T_1(h)T_2(h)\cdots T_N(h).
\end{equation}
At the state level, this corresponds to applying the frozen maps in the chronological order $\Phi_1(h,\cdot), \Phi_2(h,\cdot), \ldots, \Phi_N(h,\cdot)$. The product formula is
\begin{equation}
\label{LTH}
T(t)g
=
\tau\text{-}\lim_{n\to\infty}
\left(
T_1\!\left(\frac{t}{n}\right)
T_2\!\left(\frac{t}{n}\right)
\cdots
T_N\!\left(\frac{t}{n}\right)
\right)^n g .
\end{equation}
For $N=2$, this is
\begin{equation}
\label{LT}
T(t)g
=
\tau\text{-}\lim_{n\to\infty}
\left(
T_1\!\left(\frac{t}{n}\right)
T_2\!\left(\frac{t}{n}\right)
\right)^n g .
\end{equation}
Under the usual regularity assumptions for first-order splitting, the global error for fixed final time is of order $O(n^{-1})$, or equivalently $O(h)$ with $h=t/n$.

The Strang one-step operator is the symmetric product
\begin{equation}
\label{SPHstep}
V_{\mathrm{S}}(h)
=
\left(\prod_{i=1}^{N-1} T_i\!\left(\frac{h}{2}\right)\right)
T_N(h)
\left(\prod_{i=N-1}^{1} T_i\!\left(\frac{h}{2}\right)\right),
\end{equation}
where products are ordered from left to right in the index order.

Thus, the state updates are applied in the chronological order $1,2,\ldots,N-1,N,N-1,\ldots,2,1$, with half steps except for the middle $N$-th coordinate step. The corresponding product formula is
\begin{equation}
\label{SPH}
T(t)g
=
\tau\text{-}\lim_{n\to\infty}
V_{\mathrm{S}}\!\left(\frac{t}{n}\right)^n g .
\end{equation}
For $N=2$, this reduces to
\begin{equation}
\label{SP}
T(t)g
=
\tau\text{-}\lim_{n\to\infty}
\left(
T_1\!\left(\frac{t}{2n}\right)
T_2\!\left(\frac{t}{n}\right)
T_1\!\left(\frac{t}{2n}\right)
\right)^n g .
\end{equation}
Under the usual regularity assumptions for second-order splitting, the global error for fixed final time is of order $O(n^{-2})$, or equivalently $O(h^2)$ with $h=t/n$.

The ordering conventions in \eqref{LTstep} and \eqref{SPHstep} are not unique. Reversing the coordinate order gives another splitting of the same formal order. The algorithms below state the order of frozen updates used in the numerical implementation.

All formulas in this subsection involve nonnegative real time steps. Products with complex coefficients require additional assumptions, because $T_i(z)$ is not generally defined for complex $z$ for a general Koopman semigroup on $C_b(\Omega)$. We address this restriction in the next subsection.

\subsection{Higher-order exponential splitting schemes}
\label{sec:higher}

We now recall the background on the higher-order exponential compositions used in the numerical experiments. The formulas are classical splitting formulas for two operators, and the later algorithms use them as building blocks for coordinate-wise state updates.

Let $X$ be a Banach space, and suppose that $\mathcal{K}=\mathcal{K}_1+\mathcal{K}_2$ is the decomposition of a generator into two components. A splitting method with $s$ stages has the form
\begin{equation}
\label{sec:HigherOrder}
V_p(t)
=
\prod_{j=1}^s
e^{\alpha_j t\mathcal{K}_1}
e^{\beta_j t\mathcal{K}_2}.
\end{equation}
The coefficients $\alpha_j$ and $\beta_j$ are chosen so that the method has algebraic order $p$. This means that, when the operators are replaced by bounded matrices $\mathcal{M}_1$ and $\mathcal{M}_2$, with $\mathcal{M}=\mathcal{M}_1+\mathcal{M}_2$, the expansion satisfies $V_p(t)-e^{t\mathcal{M}} = \mathcal{O}(t^{p+1})$ as $t\to 0$. The algebraic order condition is independent of the particular ODE examples. Convergence for unbounded operators requires additional domain and stability assumptions.

For real nonnegative coefficients, the factors in \eqref{sec:HigherOrder} are semigroup factors. For complex coefficients, the notation requires more care. A general Koopman semigroup on $C_b(\Omega)$ is defined for real times $t\ge 0$, without in general defining $e^{z\mathcal{K}_i}$ for $z\in\mathbb{C}$. Thus, whenever a complex coefficient occurs, one must assume either an analytic semigroup setting or a problem-specific complex-time extension of the frozen state maps. In the numerical examples of this paper, complex-coefficient schemes are used only for frozen scalar maps that admit such complex-time evaluation.

The following estimate summarizes the form of the higher-order splitting result used here; see~\citet{JD5,JD6}.

\begin{theorem}[Higher-order splitting estimate~\citep{JD5,JD6}]
\label{ha}
Let $\mathcal{K}_1$, $\mathcal{K}_2$, and $\mathcal{K}=\mathcal{K}_1+\mathcal{K}_2$ be generators in a setting where all factors in \eqref{sec:HigherOrder} are defined for $0\le t\le T$. If some coefficients are complex, assume that the relevant semigroups are analytic on sectors containing the rays determined by these coefficients. Suppose that $V_p$ has algebraic order $p$ and is stable on $[0,T]$, namely $\sup_{n\in\mathbb{N}}\sup_{0\le t\le T} \left\| V_p\!\left(\frac{t}{n}\right)^n \right\| < \infty$. Let $U\subset D(\mathcal{K}^{p+1})$ be a subspace such that, for every word $E_{p+1} = \mathcal{K}_{i_1}\mathcal{K}_{i_2}\ldots \mathcal{K}_{i_{p+1}},~i_\ell\in\{1,2\}$, and every $u\in U$, the expression $E_{p+1}e^{t\mathcal{K}}u$ is defined in $X$ for $0\le t\le T$ and satisfies $\sup_{0\le t\le T} \bigl\|E_{p+1}e^{t\mathcal{K}}u\bigr\| < \infty$. Then, for every $u\in U$, there is a constant $C_{u,T}$ such that
\[
\left\|
V_p\!\left(\frac{t}{n}\right)^n u
-
e^{t\mathcal{K}}u
\right\|
\le
C_{u,T} n^{-p},
~
0\le t\le T .
\]
\end{theorem}

The theorem explains the role of algebraic order in the operator setting. In the Koopman setting of this paper and in our contributed algorithms in Section~\ref{sec:alg}, the real-time formulas are interpreted through the frozen Koopman semigroups described above. Complex-time formulas are interpreted at the level of the available frozen scalar maps.

A third-order complex composition used in this paper is obtained by composing Strang-type factors with conjugate coefficients.

\begin{proposition}[Third-order complex composition~\citep{JD6,JD7}]
\label{prop333}
Let $\alpha = \frac{1}{2} + \frac{i\sqrt{3}}{6}$.
Then the composition
\begin{equation}
\label{eq2d3o}
V_3(t)
=
e^{\overline{\alpha}\frac{t}{2}\mathcal{K}_2}
e^{\overline{\alpha}t\mathcal{K}_1}
e^{\frac{t}{2}\mathcal{K}_2}
e^{\alpha t\mathcal{K}_1}
e^{\alpha\frac{t}{2}\mathcal{K}_2}
\end{equation}
is algebraically of order $3$, where the bar denotes complex conjugation.
\end{proposition}

Proposition~\ref{prop333} is an algebraic order statement that follows from the usual formal expansion of the composition in noncommuting variables, where the coefficients are chosen so that all terms up to order $3$ agree with those of $e^{t(\mathcal{K}_1+\mathcal{K}_2)}$. Its use as a numerical method for a Koopman splitting requires the complex-time interpretation described above.

Higher-order formulas may also be obtained by recursively composing lower-order formulas. Here, we recall three such recursive constructions from~\citet{JD6}. These formulas are stated for a two-operator splitting $\mathcal{K}=\mathcal{K}_1+\mathcal{K}_2$. In this paper, they are used as two-operator building blocks and as a source of coefficient patterns for the implementations. The order statements of Propositions~\ref{splitU}-\ref{SplitZ} should therefore be read in the two-operator setting, under the same analytic or complex-time assumptions discussed above.

The order and factor-count statements of Propositions~\ref{splitU}-\ref{SplitZ} are the recursive algebraic identities established in~\citet{JD6}; the analytic or complex-time assumptions are needed only to interpret the factors as operators.

\begin{proposition}[Recursive family $U_{[k]}$; see~\citet{JD6}]
\label{splitU}
Let
$U_{[0]}(t)
=
V_2(t)
=
e^{\frac{t}{2}\mathcal{K}_2}
e^{t\mathcal{K}_1}
e^{\frac{t}{2}\mathcal{K}_2}.
$
For $1\le k\le 4$, define
\begin{equation}
\label{H53}
U_{[k]}(t)
=
U_{[k-1]}(\overline{a}_k t)\,
U_{[k-1]}(a_k t),
~
a_k
=
\frac{1}{2}
+
i\,
\frac{\sin(\pi/(k+2))}
{2+2\cos(\pi/(k+2))}.
\end{equation}
Then $U_{[k]}$ is algebraically of order $k+2$ for the two-operator splitting
$\mathcal{K}=\mathcal{K}_1+\mathcal{K}_2$. After merging adjacent factors involving the same operator, $U_{[k]}$ contains $2^{k+1}+1$ exponential factors.
\end{proposition}

For example, $U_{[1]}$ is the third-order composition
\[
U_{[1]}(t)
=
e^{\overline{a}_1\frac{t}{2}\mathcal{K}_2}
e^{\overline{a}_1 t\mathcal{K}_1}
e^{\frac{t}{2}\mathcal{K}_2}
e^{a_1 t\mathcal{K}_1}
e^{a_1\frac{t}{2}\mathcal{K}_2},
\]
where $a_1=\frac{1}{2}+\frac{i\sqrt{3}}{6}$. This is the formula written in \eqref{eq2d3o}. The next member, $U_{[2]}$, has order $4$ and contains $9$ exponential factors. The family can be continued recursively, with the number of factors doubling at each step.

\begin{proposition}[Recursive family $W_{[k]}$; see~\citet{JD6}]
\label{SplitW}
Let $W_{[0]}(t)=V_2(t)$. For $1\le k\le 3$, define
\begin{equation}
\label{HO53}
W_{[k]}(t)
=
W_{[k-1]}(a_k t)\,
W_{[k-1]}((1-2a_k)t)\,
W_{[k-1]}(a_k t),
\end{equation}
where
\[
a_k
=
\frac{e^{\pi i/(2k+1)}}
{2^{1/(2k+1)}+e^{\pi i/(2k+1)}} .
\]
Then $W_{[k]}$ is algebraically of order $2k+2$ for the two-operator splitting
$\mathcal{K}=\mathcal{K}_1+\mathcal{K}_2$. After merging adjacent factors involving the same operator, $W_{[k]}$ contains $2\cdot 3^k+1$ exponential factors.
\end{proposition}

In particular, $W_{[1]}$ is a fourth-order two-operator composition with seven exponential factors:
\[
W_{[1]}(t)
=
e^{a_1\frac{t}{2}\mathcal{K}_2}
e^{a_1t\mathcal{K}_1}
e^{(1-a_1)\frac{t}{2}\mathcal{K}_2}
e^{(1-2a_1)t\mathcal{K}_1}
e^{(1-a_1)\frac{t}{2}\mathcal{K}_2}
e^{a_1t\mathcal{K}_1}
e^{a_1\frac{t}{2}\mathcal{K}_2}.
\]
Here,
$
a_1
=
\frac{e^{\pi i/3}}
{2^{1/3}+e^{\pi i/3}}$.

\begin{proposition}[Recursive family $Z_{[k]}$; see~\citet{JD6}]
\label{SplitZ}
Let $Z_{[0]}(t)=V_2(t)$. For $1\le k\le 6$, define
\begin{equation}
\label{HO533}
Z_{[k]}(t)
=
Z_{[k-1]}(a_k t)\,
Z_{[k-1]}(\overline{a}_k t)\,
Z_{[k-1]}(\overline{a}_k t)\,
Z_{[k-1]}(a_k t),
\end{equation}
where
\[
a_k
=
\frac{1}{4}
+
i\,
\frac{\sin(\pi/(2k+1))}
{4+4\cos(\pi/(2k+1))}.
\]
Then $Z_{[k]}$ is algebraically of order $2k+2$ for the two-operator splitting
$\mathcal{K}=\mathcal{K}_1+\mathcal{K}_2$. After merging adjacent factors involving the same operator, $Z_{[k]}$ contains $2\cdot 4^k+1$ exponential factors.
\end{proposition}

For example, $Z_{[1]}$ is a fourth-order two-operator composition with nine exponential factors:
\[
Z_{[1]}(t)
e^{a_1\frac{t}{2}\mathcal{K}_2}
e^{a_1t\mathcal{K}_1}
e^{\frac{t}{4}\mathcal{K}_2}
e^{\overline{a}_1t\mathcal{K}_1}
e^{\overline{a}_1t\mathcal{K}_2}
e^{\overline{a}_1t\mathcal{K}_1}
e^{\frac{t}{4}\mathcal{K}_2}
e^{a_1t\mathcal{K}_1}
e^{a_1\frac{t}{2}\mathcal{K}_2}.
\]
Here, $a_1=\frac{1}{4}+\frac{i\sqrt{3}}{12}$.

The recursive definitions above are sufficient to generate the higher-order two-operator compositions used in the numerical implementation. They also determine the corresponding number of exponential factors in two dimensions. For systems with more than two coordinate generators, the next subsection provides the factor counts produced by the recursive substitution procedure used in the implementation. Those counts are reported as costs of that construction.

\subsection{Factor counts for higher-dimensional coordinate compositions}
\label{expHI}

We next state the number of elementary frozen-flow factors produced by the recursive coordinate construction used in the implementation. This subsection provides the relevant bookkeeping.

For Lie-Trotter splitting in dimension $N$, one time step contains $N$ frozen-flow factors. For Strang splitting, one time step contains $2N-1$ frozen-flow factors. These counts follow directly from the formulas \eqref{LTstep} and \eqref{SPHstep}.

For higher-order formulas, the implementation starts from a two-operator composition and recursively replaces one coordinate factor by an alternating block involving a new coordinate. We describe only the factor count. Let $q$ be the number of factors in the two-operator base formula. In the examples used here, $q$ is odd. Thus the two operators occur $\frac{q-1}{2}$ and $\quad\frac{q+1}{2}$ times, respectively, after adjacent factors of the same type have been merged.

To pass from $m-1$ coordinates to $m$ coordinates, choose a coordinate factor with the smallest current multiplicity. If that factor occurs $r$ times, replace each of its occurrences by the same alternating two-operator pattern, now using the old coordinate and the new coordinate. At the level of counts, this removes $r$ occurrences of the old factor and replaces them by $r\frac{q-1}{2}$ old-coordinate factors and $r\frac{q+1}{2}$ new-coordinate factors. This gives Algorithm~\ref{alg:expterms}.

\begin{algorithm}[t!]
\caption{Factor count for the recursive coordinate substitution}
\label{alg:expterms}
\begin{algorithmic}[1]
\STATE \textbf{Input:} Target dimension $N\ge 2$ and odd two-operator factor count $q$
\STATE Set $a=(q-1)/2$ and $b=(q+1)/2$
\STATE Initialize the list of multiplicities $L=[a,b]$
\FOR{$m=3$ to $N$}
    \STATE Sort $L$ in nondecreasing order
    \STATE Remove one smallest entry $r$ from $L$
    \STATE Append $ar$ and $br$ to $L$
\ENDFOR
\STATE \textbf{Output:} $\sum_{\ell\in L}\ell$
\end{algorithmic}
\end{algorithm}

For $N=2$, the loop in Algorithm~\ref{alg:expterms} is empty, and the output is $a+b=q$, the factor count of the base two-operator formula. For example, a two-operator third-order composition with $q=5$ factors has multiplicities $2$ and $3$. The recursive count therefore gives $5,13,25,41,65,\ldots$ factors in dimensions $2,3,4,5,6,\ldots$, respectively. Similarly, starting from a two-operator fourth-order composition with $q=7$ factors gives $7,25,49,103,175,\ldots$.

The recursive substitution rule preserves the nested coordinate structure of the implementation. Table~\ref{tabNterms} provides the number of frozen-flow factors generated by this construction. The counts are not claimed to be minimal.

\begin{table}[t!]
\centering
\caption{Number of elementary frozen-flow factors produced by the recursive substitution rule. Empty entries correspond to cases not used in the numerical experiments.}
\label{tabNterms}
\setlength{\tabcolsep}{4pt}
\fontsize{8}{8}\selectfont
\begin{tabular}{lrrrrrrrrr}
\toprule
\multicolumn{1}{c}{\textbf{Base formula}} &
\multicolumn{1}{c}{$\mathbf{p}$} &
\multicolumn{1}{c}{$\mathbf{N=2}$} &
\multicolumn{1}{c}{$\mathbf{N=3}$} &
\multicolumn{1}{c}{$\mathbf{N=4}$} &
\multicolumn{1}{c}{$\mathbf{N=5}$} &
\multicolumn{1}{c}{$\mathbf{N=6}$} &
\multicolumn{1}{c}{$\mathbf{N=7}$} &
\multicolumn{1}{c}{$\mathbf{N=8}$}\\
\midrule
\rowcolor{gray!10} Lie-Trotter & $1$  & $2$ & $3$ & $4$ & $5$ & $6$ & $7$ & $8$ \\
\rowcolor{white} Strang          & $2$  & $3$ & $5$ & $7$ & $9$ & $11$ & $13$ & $15$ \\
\rowcolor{gray!10} $3$-rd order   & $3$  & $5$ & $13$ & $25$ & $41$ & $65$ & $89$ & $121$ \\
\rowcolor{white} $4$-th order    & $4$  & $7$ & $25$ & $49$ & $103$ & $175$ & $247$ & $343$ \\
\rowcolor{gray!10} $6$-th order   & $6$  & $17$ & $145$ & $289$ & $1313$ & $2465$ & $3617$ & $4913$ \\
\rowcolor{white} $8$-th order    & $8$  & $55$ & $1513$ & $3025$ & $42391$ & $83215$ & $124039$ & $166375$ \\
\rowcolor{gray!10} $10$-th order   & $10$ & $513$ & $131585$ & $263169$ & $33817601$ & $67503105$ & & \\
\rowcolor{white} $12$-th order   & $12$ & $2049$ & $2099201$ & $4198401$ & & & & \\
\rowcolor{gray!10} $14$-th order & $14$ & $8193$ & $33562625$ & & & & & \\
\bottomrule
\end{tabular}
\end{table}

The rapid growth in Table~\ref{tabNterms} explains the range of dimensions and orders used in the experiments. In two dimensions, high-order compositions remain feasible for the examples considered here. In three dimensions, the number of factors grows much faster, and the numerical experiments are therefore restricted to lower orders. Table~\ref{tabNterms} quantifies this implementation cost.

\subsection{Algorithms}
\label{sec:alg}

We now write the coordinate-wise splitting methods as state-update algorithms. The Koopman notation identifies the frozen semigroups being composed, but the actual computation advances points in the state space. Thus, the output of the algorithm is an approximate state value. Observables are evaluated afterward by composition with this approximate state map.

\subsubsection{Two-dimensional initial value problems}

Let $O\subseteq\mathbb{R}^2$ be open, let $\Omega\subseteq O$, and let $F=(F_1,F_2):O\to\mathbb{R}^2$. Consider the initial-value problem
\begin{equation}
\label{IVP2}
x'(t)=F_1(x(t),y(t)),
~
y'(t)=F_2(x(t),y(t)),
~
(x(0),y(0))=(x,y)\in\Omega .
\end{equation}
Assume that \eqref{IVP2} generates a unique global, jointly continuous, $\Omega$-invariant flow $\varphi_t(x,y)=(x(t),y(t))$. The associated Koopman semigroup is
\[
T(t)g(x,y)=g(\varphi_t(x,y)),
~g\in C_b(\Omega).
\]

We define the frozen scalar flows as follows. For fixed $y$, let $\sigma_y(s,x)$ be the solution of
\[
\frac{d}{ds}X(s)=F_1(X(s),y),
~
X(0)=x.
\]
For fixed $x$, let $\gamma_x(s,y)$ be the solution of
\[
\frac{d}{ds}Y(s)=F_2(x,Y(s)),
~
Y(0)=y.
\]
We assume that these frozen flows are defined for the times used here and that their state maps preserve $\Omega$. The corresponding frozen state maps are
\[
\Phi_1(s,(x,y))=(\sigma_y(s,x),y),
~
\Phi_2(s,(x,y))=(x,\gamma_x(s,y)).
\]
Therefore, for $g\in C_b(\Omega)$,
\[
T_1(s)g(x,y)
=
g(\sigma_y(s,x),y),
~
T_2(s)g(x,y)
=
g(x,\gamma_x(s,y)).
\]
This notation keeps the frozen coordinate in the argument of $g$.

Let $t>0$ be the final time and let $n\in\mathbb{N}$. We write $h=\frac{t}{n}$ for the uniform time step. The Lie-Trotter state update is $\Psi_{\mathrm{LT},h} = \Phi_2(h,\cdot)\circ\Phi_1(h,\cdot)$. Equivalently, starting from $(x_0,y_0)=(x,y)$, for $k=0,\ldots,n-1$,
\begin{center}
\begin{minipage}{0.44\textwidth}
\[
\begin{alignedat}{3}
{{\color{red}(1)}}&\quad x_{k+1}
&{}={}&
\sigma_{y_k}(h,x_k),\\
{{\color{blue}(2)}}&\quad y_{k+1}
&{}={}&
\gamma_{x_{k+1}}(h,y_k).
\end{alignedat}
\]
\end{minipage}\hfill
\begin{minipage}{0.54\textwidth}
    \includegraphics[keepaspectratio, scale=.7]{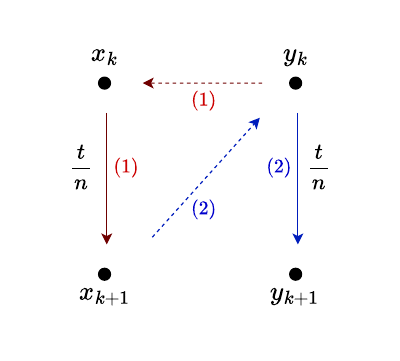}
\end{minipage}
\end{center}
The corresponding observable product is $V_{\mathrm{LT}}(h)=T_1(h)T_2(h)$, and
\[
V_{\mathrm{LT}}(h)g(x_k,y_k)
=
g(x_{k+1},y_{k+1}).
\]
For first-order splitting, the corresponding error estimate gives a global error of order $O(h)$ at fixed final time, or equivalently $O(n^{-1})$.

The Strang state update is
\[
\Psi_{\mathrm{S},h}
=
\Phi_1\!\left(\frac{h}{2},\cdot\right)
\circ
\Phi_2(h,\cdot)
\circ
\Phi_1\!\left(\frac{h}{2},\cdot\right).
\]
Thus, starting from $(x_0,y_0)=(x,y)$, for $k=0,\ldots,n-1$,
\begin{center}
\begin{minipage}{0.44\textwidth}
\[
\begin{alignedat}{3}
{{\color{red}(1)}}&\quad x_{k+\frac12}
&{}={}&
\sigma_{y_k}\!\left(\frac{h}{2},x_k\right),\\
{{\color{blue}(2)}}&\quad y_{k+1}
&{}={}&
\gamma_{x_{k+\frac12}}(h,y_k),\\
{{\color{green}(3)}}&\quad x_{k+1}
&{}={}&
\sigma_{y_{k+1}}\!\left(\frac{h}{2},x_{k+\frac12}\right).
\end{alignedat}
\]
\end{minipage}\hfill
\begin{minipage}{0.54\textwidth}
    \includegraphics[keepaspectratio, scale=.7]{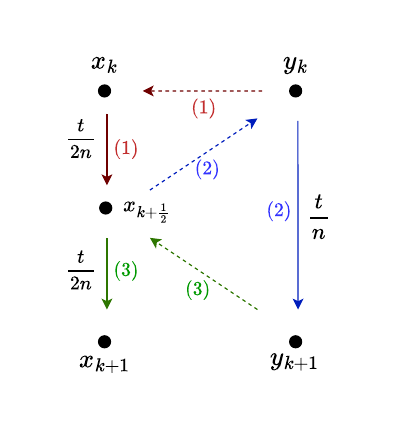}
\end{minipage}
\end{center}
The corresponding observable product is
\[
V_{\mathrm{S}}(h)
=
T_1\!\left(\frac{h}{2}\right)
T_2(h)
T_1\!\left(\frac{h}{2}\right),
\]
and
\[
V_{\mathrm{S}}(h)g(x_k,y_k)
=
g(x_{k+1},y_{k+1}).
\]
For second-order splitting, the corresponding error estimate gives a global error of order $O(h^2)$ at fixed final time, or equivalently $O(n^{-2})$.

Higher-order two-dimensional schemes are implemented by composing the same frozen maps with the coefficients described in Section~\ref{sec:higher}. The third-order two-dimensional update used as a template for the higher-order implementations is given in~\ref{supp:supplementary_algo}.

If a coefficient is complex, the corresponding step is used only when the frozen maps $\sigma_y$ and $\gamma_x$ admit the required complex-time evaluation. In such cases, intermediate states may be complex. The final reported state is interpreted according to the problem-specific implementation described in the numerical section.

All updates above assume exact evaluation of the frozen scalar maps. If the scalar frozen problems are solved numerically, then the scalar-solver error must be added to the splitting error. In the numerical examples in Section~\ref{sec:experiments}, the frozen maps used in the computations are available in closed form.

\subsubsection{Three-dimensional initial value problems}

Let $O\subseteq\mathbb{R}^3$ be open, let $\Omega\subseteq O$, and let $F=(F_1,F_2,F_3):O\to\mathbb{R}^3$. Consider the initial-value problem
\begin{align*}
x'(t)&=F_1(x(t),y(t),z(t)),\\
y'(t)&=F_2(x(t),y(t),z(t)),\\
z'(t)&=F_3(x(t),y(t),z(t)),
\end{align*}
with initial condition $(x(0),y(0),z(0))=(x,y,z)\in\Omega$. We assume that this problem generates a unique global, jointly continuous, $\Omega$-invariant flow $\varphi_t(x,y,z)=(x(t),y(t),z(t))$. The associated Koopman semigroup is
\[
T(t)g(x,y,z)=g(\varphi_t(x,y,z)),
~
g\in C_b(\Omega).
\]

We define the frozen scalar flows by fixing all variables except one. For fixed $y$ and $z$, let $\sigma_{y,z}(s,x)$ solve
\[
\frac{d}{ds}X(s)=F_1(X(s),y,z),
~
X(0)=x.
\]
For fixed $x$ and $z$, let $\gamma_{x,z}(s,y)$ solve
\[
\frac{d}{ds}Y(s)=F_2(x,Y(s),z),
~
Y(0)=y.
\]
For fixed $x$ and $y$, let $\tau_{x,y}(s,z)$ solve
\[
\frac{d}{ds}Z(s)=F_3(x,y,Z(s)),
~
Z(0)=z.
\]
We assume that these frozen flows are defined for the times used below and that their state maps preserve $\Omega$. The corresponding frozen state maps are
\begin{align*}
\Phi_1(s,(x,y,z))&=(\sigma_{y,z}(s,x),y,z),\\
\Phi_2(s,(x,y,z))&=(x,\gamma_{x,z}(s,y),z),\\
\Phi_3(s,(x,y,z))&=(x,y,\tau_{x,y}(s,z)).
\end{align*}
Thus, for $g\in C_b(\Omega)$,
\begin{align*}
T_1(s)g(x,y,z)
&=
g(\sigma_{y,z}(s,x),y,z),\\
T_2(s)g(x,y,z)
&=
g(x,\gamma_{x,z}(s,y),z),\\
T_3(s)g(x,y,z)
&=
g(x,y,\tau_{x,y}(s,z)).
\end{align*}
This notation keeps the frozen coordinates in the argument of $g$.

Let $t>0$ be the final time, let $n\in\mathbb{N}$, and set $h=t/n$. We use the coordinate order $z,y,x$ in the Lie-Trotter update. The state update is
\[
\Psi_{\mathrm{LT},h}^{(3)}
=
\Phi_1(h,\cdot)\circ\Phi_2(h,\cdot)\circ\Phi_3(h,\cdot).
\]
Thus, starting from $(x_0,y_0,z_0)=(x,y,z)$, for $k=0,\ldots,n-1$,
\begin{center}
\begin{minipage}{0.44\textwidth}
\[
\begin{alignedat}{3}
{{\color{red}(1)}}&\quad z_{k+1}
&{}={}&
\tau_{x_k,y_k}(h,z_k),\\
{{\color{blue}(2)}}&\quad y_{k+1}
&{}={}&
\gamma_{x_k,z_{k+1}}(h,y_k),\\
{{\color{green}(3)}}&\quad x_{k+1}
&{}={}&
\sigma_{y_{k+1},z_{k+1}}(h,x_k).
\end{alignedat}
\]
\end{minipage}
\begin{minipage}{0.54\textwidth}
    \begin{center}
    \includegraphics[keepaspectratio, scale=1.2]{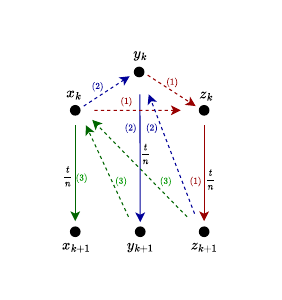}
    \end{center}
\end{minipage}
\end{center}
The corresponding observable product is
\[
V_{\mathrm{LT}}^{(3)}(h)
=
T_3(h)T_2(h)T_1(h),
\]
and
\[
V_{\mathrm{LT}}^{(3)}(h)g(x_k,y_k,z_k)
=
g(x_{k+1},y_{k+1},z_{k+1}).
\]
For first-order splitting, the corresponding error estimate gives a global error of order $O(h)$ at fixed final time, equivalently $O(n^{-1})$.

For completeness, we also state the corresponding three-dimensional Strang update in the same coordinate order. The state update is
\[
\Psi_{\mathrm{S},h}^{(3)}
=
\Phi_3\!\left(\frac{h}{2},\cdot\right)
\circ
\Phi_2\!\left(\frac{h}{2},\cdot\right)
\circ
\Phi_1(h,\cdot)
\circ
\Phi_2\!\left(\frac{h}{2},\cdot\right)
\circ
\Phi_3\!\left(\frac{h}{2},\cdot\right).
\]
Equivalently, for $k=0,\ldots,n-1$,
\[
\begin{aligned}
z_{k+\frac12}
&=
\tau_{x_k,y_k}\!\left(\frac{h}{2},z_k\right),\\
y_{k+\frac12}
&=
\gamma_{x_k,z_{k+\frac12}}\!\left(\frac{h}{2},y_k\right),\\
x_{k+1}
&=
\sigma_{y_{k+\frac12},z_{k+\frac12}}(h,x_k),\\
y_{k+1}
&=
\gamma_{x_{k+1},z_{k+\frac12}}\!\left(\frac{h}{2},y_{k+\frac12}\right),\\
z_{k+1}
&=
\tau_{x_{k+1},y_{k+1}}\!\left(\frac{h}{2},z_{k+\frac12}\right).
\end{aligned}
\]
The corresponding product is
\[
V_{\mathrm{S}}^{(3)}(h)
=
T_3\!\left(\frac{h}{2}\right)
T_2\!\left(\frac{h}{2}\right)
T_1(h)
T_2\!\left(\frac{h}{2}\right)
T_3\!\left(\frac{h}{2}\right).
\]
For second-order splitting, the corresponding error estimate gives a global error
of order $O(h^2)$ at fixed final time, equivalently $O(n^{-2})$.

Higher-order three-dimensional schemes are obtained by composing the frozen maps $\Phi_1$, $\Phi_2$, and $\Phi_3$ with the coefficient expressions described in Subsection~\ref{sec:higher} and with the recursive coordinate substitution counted in Subsection~\ref{expHI}. The third-order update used in the implementation is given in Appendix~\ref{supp:supplementary_algo}. The factor counts in Table~\ref{tabNterms} describe the number of frozen-flow evaluations required by that construction.

The same restrictions for the two-dimensional case apply here. If all coefficients are real and nonnegative, the frozen maps are real-time flow maps. If complex coefficients are used, then the frozen scalar maps must admit the required complex-time evaluations. 


%% file: experiments.tex
\section{Numerical experiments}
\label{sec:experiments}

We carry out numerical experiments for three finite-dimensional ODE systems. The purpose is to evaluate the coordinate-wise splitting algorithms described above when the frozen subflows are available in closed form. The examples are the Lotka-Volterra system, the Van der Pol oscillator, and the Lorenz system.

The experiments are model-known. The vector field is given, and the splitting methods are applied to the corresponding frozen coordinate flows. No data-driven approximation of a Koopman operator or Koopman generator is constructed. The Koopman notation is used only to organize the product formulas that lead to the state-space algorithms.

For each experiment, a high-accuracy RK45 computation is used as a reference solution. This reference is not treated as an exact solution. It is used only as a common benchmark for comparing the splitting schemes on the same output grid. If $u_j^{\mathrm{split}}$ denotes the splitting approximation and $u_j^{\mathrm{ref}}$ denotes the RK45 reference value at the same time $t_j$, then the root mean square error is
\[
\mathrm{RMSE}
=
\left(
\frac{1}{m}
\sum_{j=1}^m
\|u_j^{\mathrm{split}}-u_j^{\mathrm{ref}}\|_2^2
\right)^{1/2}.
\]
In all tables, the RMSE is computed at the $m=n$ positive uniform grid points $t_j=\frac{jT}{n},~j=1,\ldots,n$, where $T$ is the final time and $n$ is the number of splitting steps.

Errors in the range $10^{-11}$-$10^{-13}$ are close to the attainable precision of the computation. In this range, differences between methods should be interpreted cautiously, since they may reflect numerical roundoff.

\subsection{Lotka-Volterra system}
\label{subsec:bio}

We first consider the Lotka-Volterra predator-prey system
\begin{align*}
x'(t)&=\alpha x(t)-\beta x(t)y(t),\\
y'(t)&=\delta x(t)y(t)-\tau y(t),
\end{align*}
where $x(0)=x_0,~y(0)=y_0$. Here, $x(t)$ denotes the prey population and $y(t)$ denotes the predator population. The vector field is
\[
F(x,y)=
\left(
x(\alpha-\beta y),
y(\delta x-\tau)
\right).
\]

The coordinate-wise frozen equations are solvable. If $y$ is held fixed, then $X'(s)=X(s)(\alpha-\beta y), ~ X(0)=x$, hence $\sigma_y(s,x) = x\exp (s(\alpha-\beta y))$. If $x$ is held fixed, then $Y'(s)=Y(s)(\delta x-\tau), ~ Y(0)=y$, hence $\gamma_x(s,y) = y\exp (s(\delta x-\tau))$. These formulas define the frozen maps used by the splitting algorithms. Since they are functions of the time parameter $s$, they also provide the complex-time evaluations needed for the complex-coefficient schemes in this example.

\begin{figure}[t!]
\centering
\includegraphics[width=8cm]{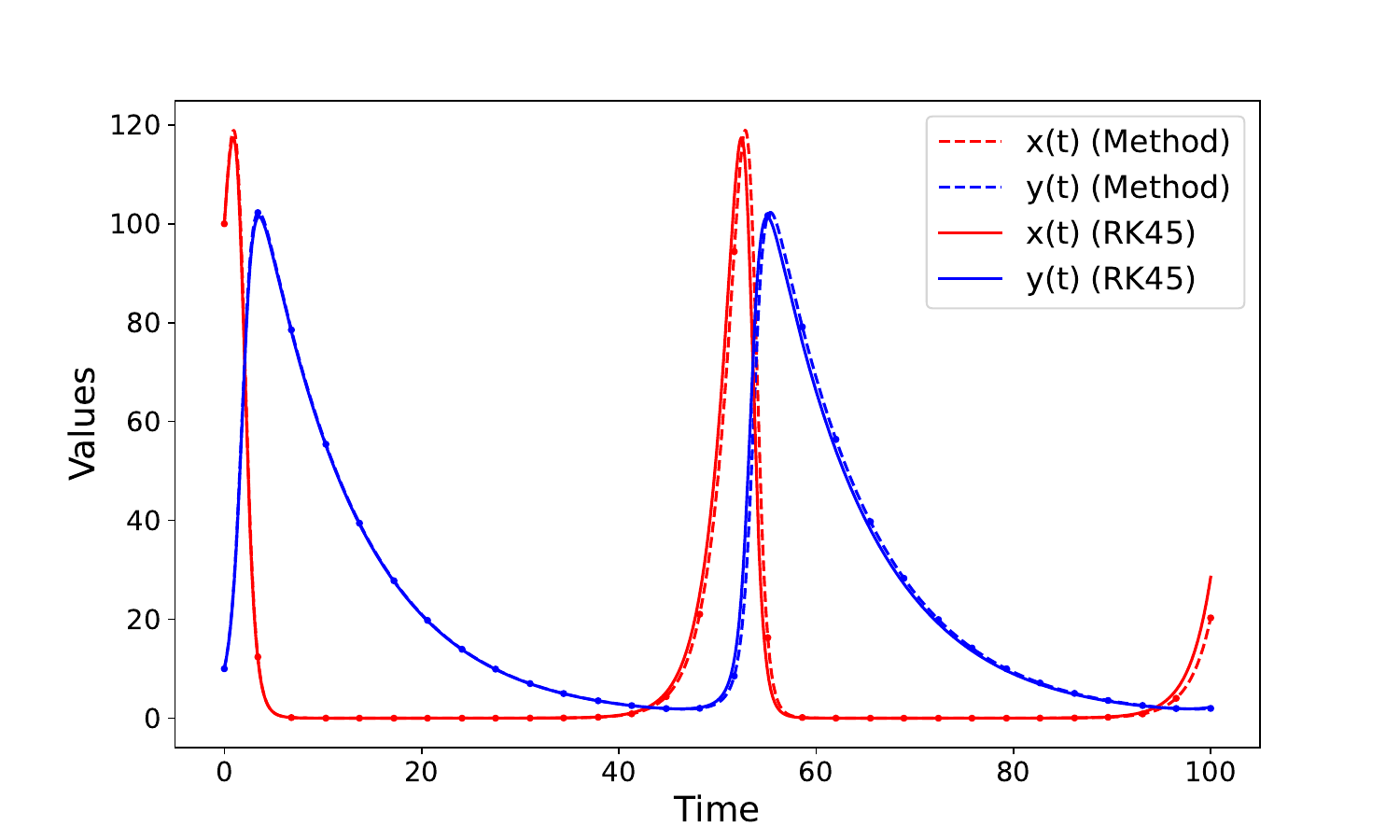}
\caption{Lotka-Volterra trajectories computed by Lie-Trotter splitting and by the RK45 reference solution over $[0,100]$ with $n=1000$ uniform splitting steps.}
\label{lok1}
\end{figure}

The numerical experiment uses $\alpha=0.5, ~ \beta=0.02, ~ \delta=0.01, ~ \tau=0.1$, with initial condition $x_0=100, ~ y_0=10$. The integration interval is $[0,100]$. Figure~\ref{lok1} shows the Lie-Trotter approximation with $n=1000$ uniform time steps, together with the RK45 reference solution evaluated on the same grid.

The RMSE values for the Lotka-Volterra experiment are reported in Table~\ref{lokRmse}. The table compares Lie-Trotter, Strang, and higher-order splitting schemes for several choices of $n$. The higher-order methods use the same frozen maps $\sigma_y$ and $\gamma_x$, composed with the coefficients described in Subsection~\ref{sec:higher}.

\begin{table}[t!]
\centering
\caption{RMSE between the Lotka-Volterra splitting approximations and the RK45 reference solution, evaluated at the same $m=n$ positive uniform grid points.}
\label{lokRmse}
\setlength{\tabcolsep}{4pt}
\fontsize{7}{8}\selectfont
\begin{tabular}{rrrrrrrrr}
\toprule
\multicolumn{1}{c}{$\mathbf{m=n}$} &
\multicolumn{1}{c}{\textbf{Lie-Trotter}} &
\multicolumn{1}{c}{\textbf{Strang}} &
\multicolumn{1}{c}{\textbf{3rd}} &
\multicolumn{1}{c}{\textbf{6th}} &
\multicolumn{1}{c}{\textbf{8th}} &
\multicolumn{1}{c}{\textbf{10th}} &
\multicolumn{1}{c}{\textbf{12th}} &
\multicolumn{1}{c}{\textbf{14th}}\\
\midrule
\rowcolor{gray!10}
$100$
& $3.30\times 10^{1}$
& $1.88\times 10^{0}$
& $1.47\times 10^{0}$
& $8.00\times 10^{-4}$
& $7.00\times 10^{-8}$
& $5.77\times 10^{-12}$
& $7.75\times 10^{-13}$
& $1.01\times 10^{-11}$ \\
\rowcolor{white}
$1000$
& $4.16\times 10^{0}$
& $1.74\times 10^{-2}$
& $2.73\times 10^{-5}$
& $2.33\times 10^{-11}$
& $1.57\times 10^{-11}$
& $1.43\times 10^{-11}$
& $9.32\times 10^{-10}$
& $1.41\times 10^{-9}$ \\
\rowcolor{gray!10}
$10000$
& $4.15\times 10^{-1}$
& $2.00\times 10^{-4}$
& $2.02\times 10^{-9}$
& $3.28\times 10^{-11}$
& $7.90\times 10^{-11}$
& $1.32\times 10^{-9}$
& $1.34\times 10^{-9}$
& $1.29\times 10^{-8}$ \\
\bottomrule
\end{tabular}
\end{table}

Table~\ref{lokRmse} shows the expected qualitative behavior. For the lower-order methods, increasing $n$ reduces the RMSE. For small $n$, the higher-order compositions reduce the error more rapidly than Lie-Trotter or Strang splitting. This is the regime in which the higher-order formulas are most useful.

The entries with very small magnitude should not be overinterpreted. At this level, the observed differences depend on roundoff error, the RK45 reference computation, and the details of complex arithmetic in the splitting implementation. Thus, the table supports the conclusion that higher-order compositions can reach the limit of machine precision with relatively few time steps.

\begin{figure}[t!]
\centering
\includegraphics[width=9cm]{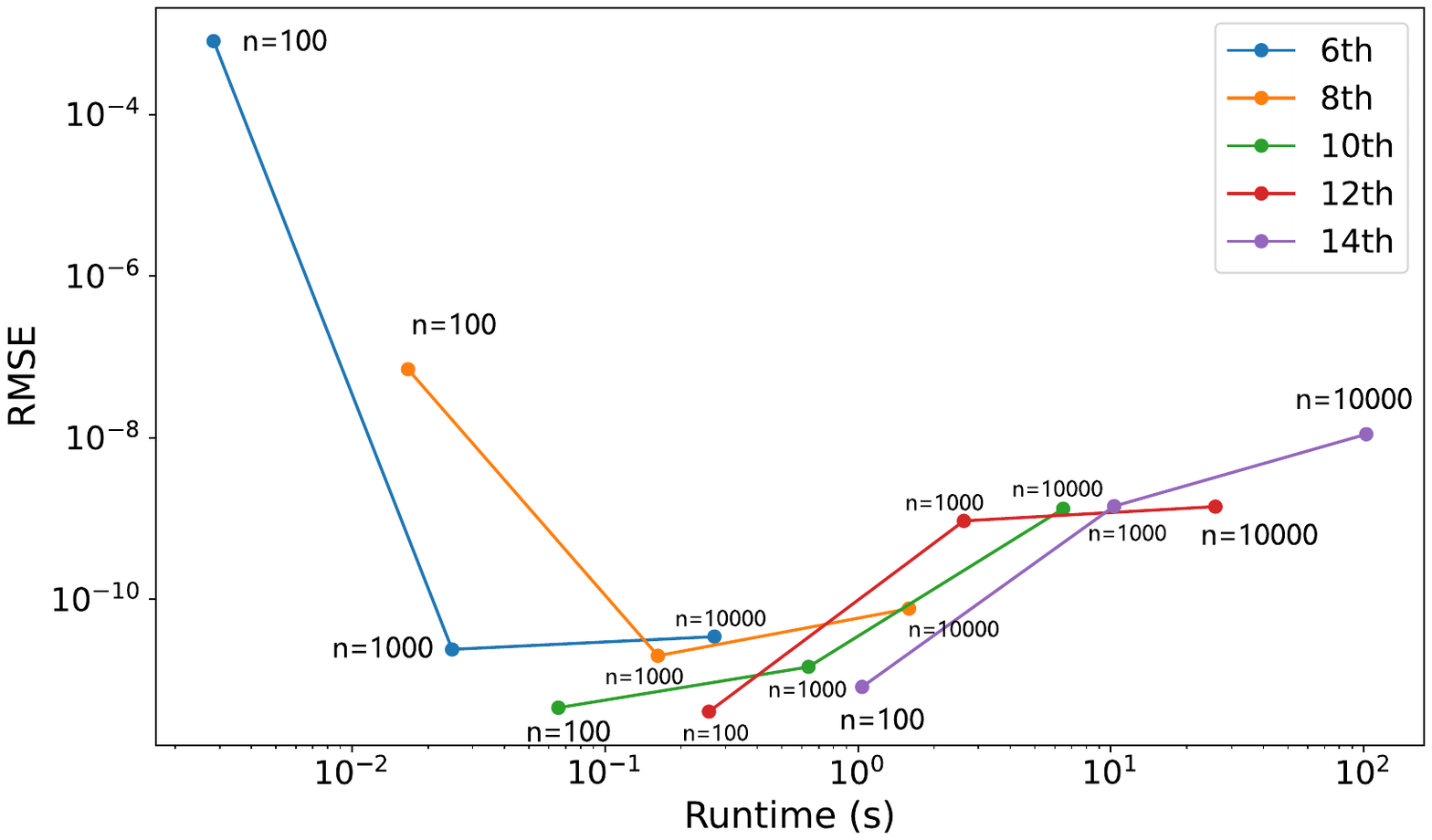}
\caption{Work-precision comparison for the Lotka-Volterra system. The RMSE relative to the RK45 reference solution is plotted against runtime in seconds.}
\label{wp1}
\end{figure}

The work-precision plot (Figure~\ref{wp1}) gives the same conclusion as Table~\ref{lokRmse} in terms of computational cost. Higher-order compositions can be advantageous when a small number of time steps is used and the frozen maps are cheap to evaluate. However, once the computation reaches the numerical precision limit, additional order does not necessarily improve the observed error.

\subsection{Van der Pol oscillator}
\label{subsec:vdp}

We next consider the Van der Pol oscillator $x''(t)+(x(t)^2-1)x'(t)+x(t)=0$. With $y(t)=x'(t)$, this equation becomes the first-order system
\begin{align*}
x'(t)&=y(t),\\
y'(t)&=(1-x(t)^2)y(t)-x(t).
\end{align*}
The initial condition used in the experiment is $x(0)=-0.2,~y(0)=0$, and the integration interval is $[0,25]$. The vector field is
\[
F(x,y)=
\left(y,
(1-x^2)y-x
\right).
\]

The frozen coordinate equations are again solvable. If $y$ is held fixed, then $X'(s)=y, ~ X(0)=x$, and therefore $\sigma_y(s,x)=x+sy $. If $x$ is held fixed, then $Y'(s)=(1-x^2)Y(s)-x,~Y(0)=y$. Writing $a(x)=1-x^2,~b(x)=-x$, the frozen $y$-flow is
\[
\gamma_x(s,y)
=
\begin{cases}
e^{s a(x)}y+
\dfrac{e^{s a(x)}-1}{a(x)}\,b(x),
& a(x)\neq 0,\\[1.0em]
y+s\,b(x),
& a(x)=0.
\end{cases}
\]
Equivalently,
\[
\gamma_x(s,y)
=
\begin{cases}
e^{s(1-x^2)}y
-
x\dfrac{e^{s(1-x^2)}-1}{1-x^2},
& x^2\neq 1,\\[1.0em]
y-sx,
& x^2=1.
\end{cases}
\]
The second formula is the removable limiting case of the first one. These frozen maps are the elementary maps used in the splitting schemes.

The dependence on the time parameter $s$ holds for both frozen maps, after the removable singularity at $a(x)=0$ is interpreted by its limit. Thus, the complex-coefficient compositions used in this example can be evaluated by applying the same formulas in complex arithmetic. Intermediate states may be complex. The reported errors are computed from the final numerical output used by the implementation and compared with the RK45 reference solution on the same grid.

\begin{figure}[t!]
\centering
\includegraphics[width=8cm]{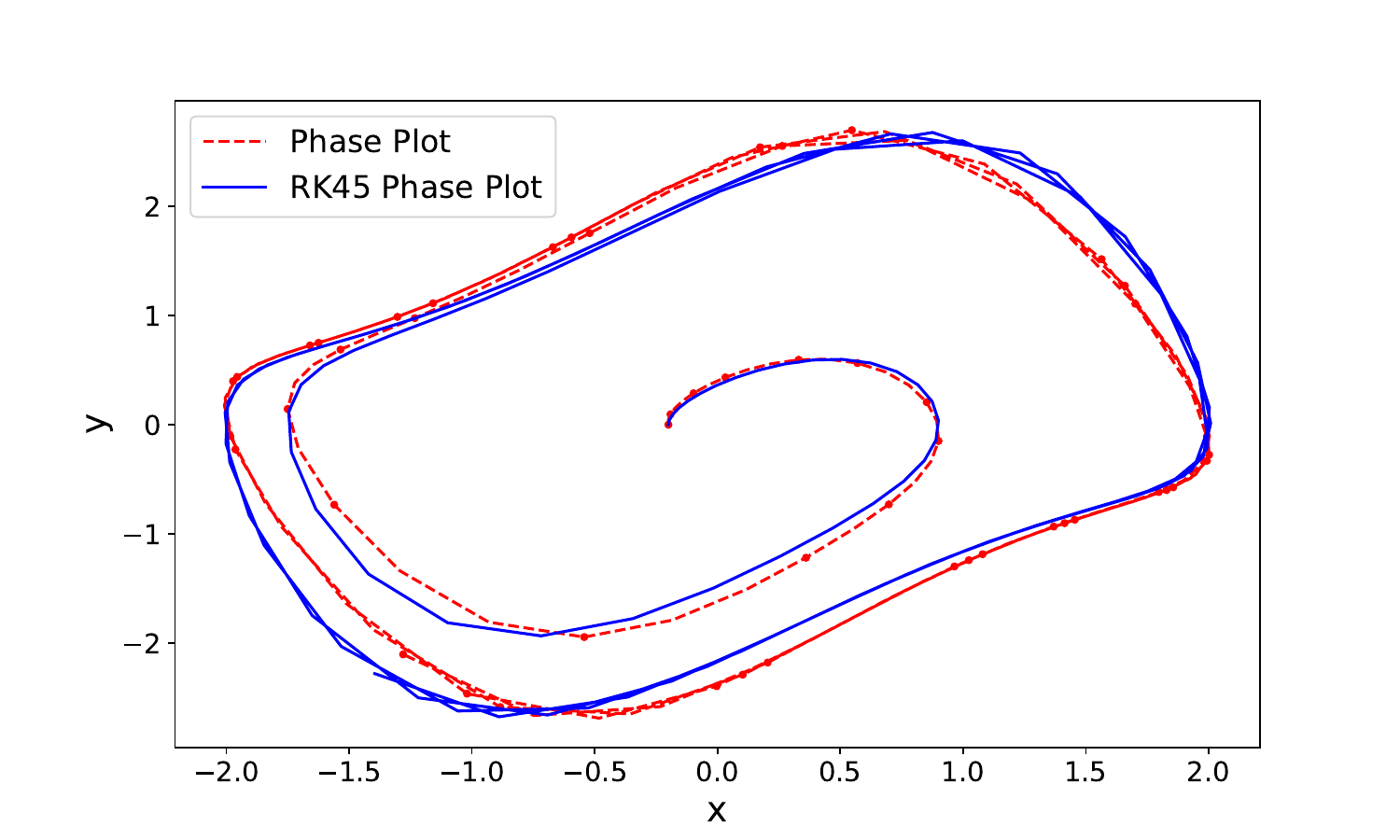}
\caption{Phase-space trajectories of the Van der Pol oscillator over $[0,25]$ with $n=125$ uniform splitting steps, compared with the RK45 reference solution evaluated on the same output grid.}
\label{van1}
\end{figure}

Figure~\ref{van1} shows the phase-space trajectory obtained with the splitting computation, together with the RK45 reference solution. The corresponding quantitative comparison is given by the RMSE values in Table~\ref{van}.

\begin{table}[t!]
\centering
\caption{RMSE between the Van der Pol splitting approximations and the RK45 reference solution, evaluated at the same $m=n$ positive uniform grid points.}
\label{van}
\setlength{\tabcolsep}{4pt}
\fontsize{7}{8}\selectfont
\begin{tabular}{rrrrrrrrr}
\toprule
\multicolumn{1}{c}{$\mathbf{m=n}$} &
\multicolumn{1}{c}{\textbf{Lie--Trotter}} &
\multicolumn{1}{c}{\textbf{Strang}} &
\multicolumn{1}{c}{\textbf{3rd}} &
\multicolumn{1}{c}{\textbf{6th}} &
\multicolumn{1}{c}{\textbf{8th}} &
\multicolumn{1}{c}{\textbf{10th}} &
\multicolumn{1}{c}{\textbf{12th}} &
\multicolumn{1}{c}{\textbf{14th}}\\
\midrule
\rowcolor{gray!10}
$125$
& $1.11\times 10^{-1}$
& $7.17\times 10^{-2}$
& $1.00\times 10^{-3}$
& $2.99\times 10^{-8}$
& $5.96\times 10^{-13}$
& $1.67\times 10^{-13}$
& $8.92\times 10^{-13}$
& $8.96\times 10^{-12}$ \\
\rowcolor{white}
$500$
& $3.10\times 10^{-2}$
& $4.30\times 10^{-3}$
& $6.84\times 10^{-6}$
& $4.70\times 10^{-12}$
& $3.64\times 10^{-13}$
& $1.50\times 10^{-12}$
& $1.11\times 10^{-11}$
& $8.49\times 10^{-11}$ \\
\rowcolor{gray!10}
$1000$
& $1.60\times 10^{-2}$
& $1.10\times 10^{-3}$
& $4.56\times 10^{-7}$
& $1.80\times 10^{-13}$
& $7.87\times 10^{-12}$
& $4.47\times 10^{-12}$
& $3.48\times 10^{-11}$
& $1.74\times 10^{-10}$ \\
\bottomrule
\end{tabular}
\end{table}

Table~\ref{van} shows the same qualitative pattern as in the Lotka-Volterra experiment. Lie-Trotter and Strang splitting improve as the number of time steps increases. The third-order method gives a substantial reduction in RMSE relative to these two lower-order schemes. The higher-order methods reach errors close to machine precision with relatively few time steps.

\begin{figure}[t!]
\centering
\includegraphics[width=9cm]{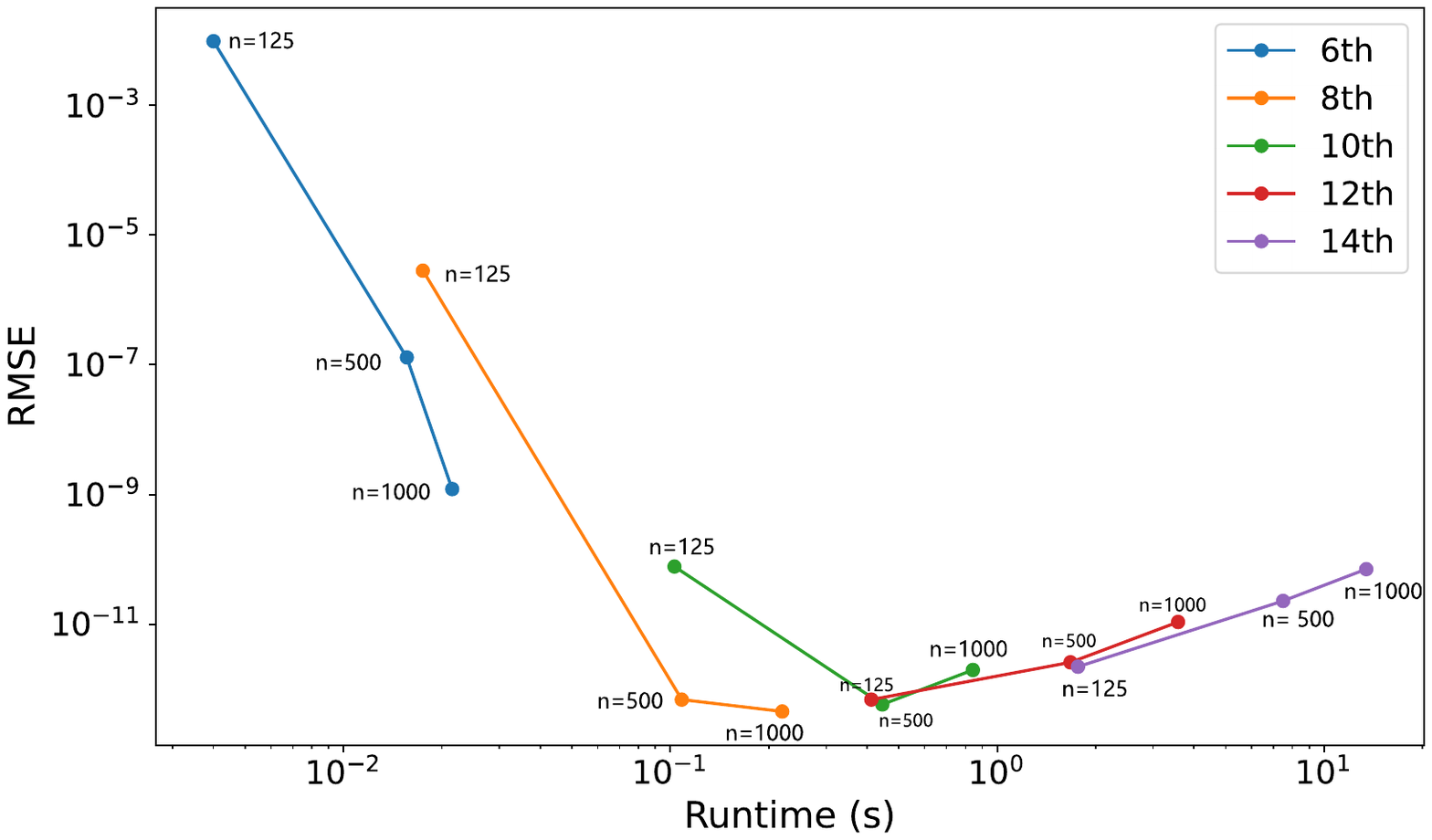}
\caption{Work-precision comparison for the Van der Pol oscillator. The RMSE relative to the RK45 reference solution is plotted against runtime in seconds.}
\label{wp2}
\end{figure}

The work-precision plot (Figure~\ref{wp2}) shows that higher-order compositions can be efficient when the frozen maps are cheap to evaluate. In this example, the high-order schemes reach small error at modest runtime. Once the errors approach the numerical precision limit, however, additional order does not necessarily produce a smaller observed RMSE.

\subsection{Lorenz system}
\label{subsec:lorenz}

We finally consider the Lorenz system
\begin{align*}
x'(t)&=\alpha(y(t)-x(t)),\\
y'(t)&=x(t)(\rho-z(t))-y(t),\\
z'(t)&=x(t)y(t)-\beta z(t).
\end{align*}
The parameters used in the experiment are $\alpha=10, ~ \rho=28, ~ \beta=\frac{8}{3}$, and the initial condition is $x(0)=1, ~ y(0)=1, ~ z(0)=1$. The integration interval is $[0,20]$.

The coordinate-wise frozen equations are affine scalar equations. If $y$ and $z$ are held fixed, then $X'(s)=\alpha(y-X(s)), ~ X(0)=x$, hence $\sigma_{y,z}(s,x) = y(1-e^{-\alpha s})+xe^{-\alpha s}$. Although the notation includes both frozen variables $y$ and $z$, this first frozen equation depends only on $y$. If $x$ and $z$ are held fixed, then $Y'(s)=x(\rho-z)-Y(s), ~ Y(0)=y$, and therefore
\[
\gamma_{x,z}(s,y) = x(\rho-z)(1-e^{-s})+ye^{-s}.
\]
Finally, if $x$ and $y$ are held fixed, then $Z'(s)=xy-\beta Z(s), ~ Z(0)=z$, so that
\[
\tau_{x,y}(s,z)
=
\frac{xy}{\beta}(1-e^{-\beta s})+ze^{-\beta s}.
\]
These three maps are the frozen state maps used in the splitting schemes.

The frozen maps are functions of the time parameter $s$. As in the previous examples, the reported errors are computed from the numerical output used in the implementation and compared with the RK45 reference solution on the same grid.

\begin{figure}[t!]
\centering
\includegraphics[width=13cm]{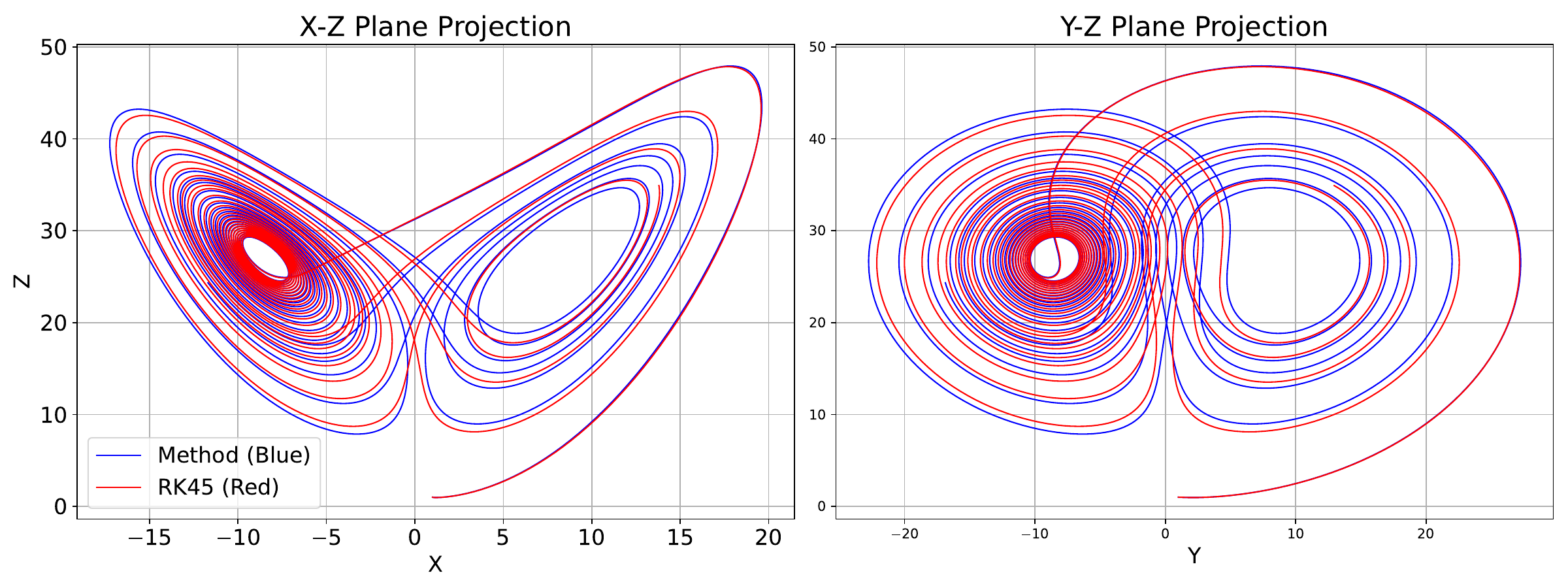}
\caption{Phase-space trajectory of the Lorenz system computed by Lie-Trotter splitting and by the RK45 reference solution over $[0,20]$ with $n=20000$ uniform splitting steps.}
\label{lor1}
\end{figure}

Figure~\ref{lor1} compares the Lie-Trotter splitting approximation with the RK45 reference solution for $n=20000$ uniform splitting steps. The figure is intended as a finite-time trajectory comparison. Since the Lorenz system is chaotic, trajectory errors over long time intervals may reflect both numerical discretization error and sensitive dependence on initial conditions. The RMSE values below should therefore be interpreted as finite-time discrepancies on the specified interval, not as long-time statistical or invariant-measure errors.

The RMSE values for the Lorenz experiment are reported in Table~\ref{lor}. The comparison uses the same definition of RMSE as in the previous examples, with the splitting approximation and the RK45 reference solution evaluated at the same $m=n$ positive uniform grid points. The table includes Lie-Trotter, Strang, third-order, and sixth-order splitting schemes. Higher orders are not included for this three-dimensional example because the recursive coordinate construction produces a rapidly increasing number of frozen-flow factors, as shown in Table~\ref{tabNterms}.

\begin{table}[t!]
\centering
\caption{RMSE between the Lorenz splitting approximations and the RK45 reference solution, evaluated at the same $m=n$ positive uniform grid points.}
\label{lor}
\setlength{\tabcolsep}{5pt}
\fontsize{9}{10}\selectfont
\begin{tabular}{rrrrr}
\toprule
\multicolumn{1}{c}{$\mathbf{m=n}$} &
\multicolumn{1}{c}{\textbf{Lie-Trotter}} &
\multicolumn{1}{c}{\textbf{Strang}} &
\multicolumn{1}{c}{\textbf{3rd}} &
\multicolumn{1}{c}{\textbf{6th}}\\
\midrule
\rowcolor{gray!10}
$1000$
& $1.549\times 10^{1}$
& $1.009\times 10^{1}$
& $7.57\times 10^{0}$
& $3.23\times 10^{-6}$ \\
\rowcolor{white}
$20000$
& $1.255\times 10^{1}$
& $1.25\times 10^{0}$
& $1.85\times 10^{-5}$
& $7.48\times 10^{-8}$ \\
\rowcolor{gray!10}
$100000$
& $1.037\times 10^{1}$
& $4.13\times 10^{-2}$
& $2.27\times 10^{-8}$
& $2.75\times 10^{-7}$ \\
\bottomrule
\end{tabular}
\end{table}

Table~\ref{lor} shows that the lower-order methods require many time steps before the finite-time RMSE becomes small. Lie-Trotter splitting remains relatively inaccurate over this interval, even for $n=100000$. Strang splitting improves substantially as $n$ increases. The third-order and sixth-order compositions give much smaller RMSE values for the tested step counts.

The sixth-order method gives the smallest RMSE for $n=1000$ and $n=20000$. For $n=100000$, the third-order method gives the smallest value in the table. This does not imply that third order is intrinsically better than sixth order for the Lorenz system. At these small error levels, the observed value can be affected by roundoff error, complex arithmetic in the higher-order composition, the RK45 reference computation, and sensitivity of the Lorenz flow. The conclusion is that higher-order coordinate-wise compositions can reduce the finite-time RMSE substantially when the frozen maps are known, but the best observed order depends on the step size and on numerical precision effects.

\begin{figure}[t!]
\centering
\includegraphics[width=9cm]{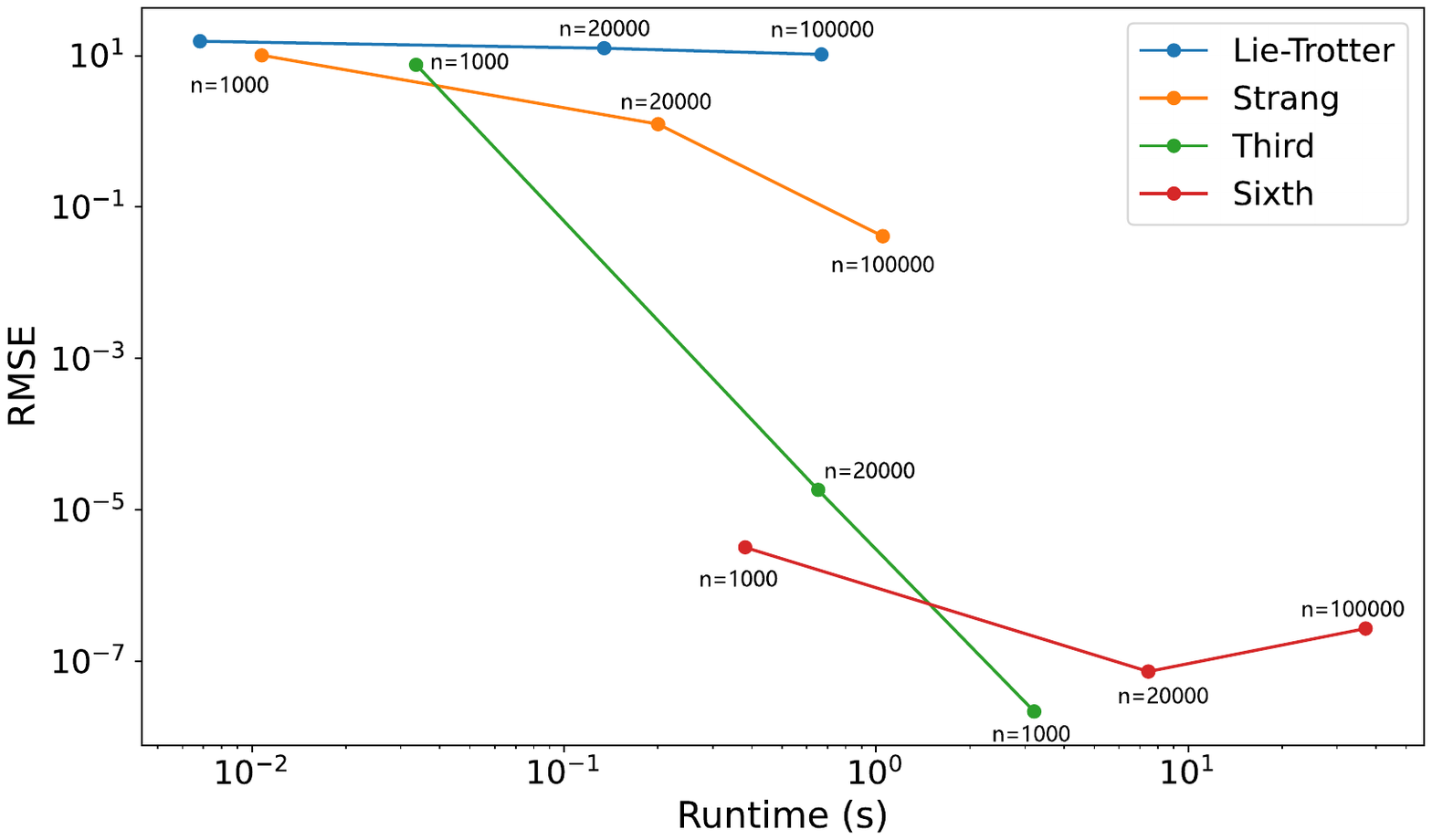}
\caption{Work-precision comparison for the Lorenz system. The RMSE relative to the RK45 reference solution is plotted against runtime in seconds. Each point corresponds to a fixed number of uniform splitting steps $n$.}
\label{wp3}
\end{figure}

Figure~\ref{wp3} reports the same comparison in terms of runtime. The plot should be read together with Table~\ref{lor}. For this implementation and these step counts, the higher-order methods reduce the RMSE more rapidly than Lie-Trotter and Strang splitting. The runtime comparison is implementation-dependent, because each method uses a different number of frozen-flow evaluations per step. In particular, the sixth-order three-dimensional composition has many more elementary factors than the lower-order schemes.

Overall, the Lorenz experiment supports the same conclusion as the two-dimensional examples. Coordinate-wise splitting with frozen maps available in closed form can give accurate finite-time approximations, and higher-order compositions can be effective when the extra frozen-flow evaluations are not too costly.

%% file: conclusion.tex
\section{Conclusion and future directions}
\label{sec:conclusions}

We have presented a coordinate-wise splitting framework for finite-dimensional ODE simulation based on Koopman-Lie product formulas and frozen scalar subflows. The main computational step is the replacement of the full flow by compositions of frozen coordinate flows. When these frozen flows are available in close form or cheap to evaluate, this gives simple, practical state-update algorithms. Lie-Trotter and Strang splitting provide the basic first- and second-order schemes. Higher-order compositions give additional accuracy when their coefficient patterns are applicable and when the corresponding frozen maps can be evaluated for the required time parameters.

The numerical experiments on the Lotka-Volterra, Van der Pol, and Lorenz systems show the finite-time behavior of these algorithms in examples where the frozen maps are available in closed form. The results illustrate the trade-off between the number of time steps, the number of frozen-flow factors, the splitting order, and runtime. In the two-dimensional examples, higher-order methods reach very small RMSE values with relatively few time steps. In the Lorenz example, higher-order compositions also reduce finite-time trajectory discrepancies on the interval considered. These numerical benchmark comparisons are made against an RK45 reference solution on the same output grid.

Future work will focus on questions that are not resolved here. One direction is to give a separate order analysis for the recursive $N$-coordinate constructions. Another is to include the error introduced when the frozen scalar subproblems are solved numerically rather than evaluated in closed form. It would also be useful to identify classes of ODEs for which the frozen subflows have real-time or complex-time formulas with favorable stability or invariant-set properties. Extensions to non-autonomous systems, local flows, and partial differential equation systems may also be possible, but these require additional assumptions on the flow maps, domains, and function spaces.

%% file: supplementary_algorithms.tex
\section{Supplementary algorithms}
\label{supp:supplementary_algo}

This appendix provides higher-order update formulas used in the implementation. The formulas are written at the level of state variables. They are obtained by applying the operator products from Subsection~\ref{sec:higher} to the frozen coordinate maps introduced in Subsection~\ref{sec:alg}.

The formulas with complex coefficients are implemented in complex arithmetic. Intermediate states may be complex. The formulas below do not imply real-valuedness, positivity, or preservation of an invariant set for the intermediate substeps.

\subsection{Two-dimensional third-order update}

We first give the third-order two-dimensional update. Let $\alpha=\frac{1}{2}+\frac{i\sqrt{3}}{6}$, let $\overline{\alpha}$ be the complex conjugate of $\alpha$, and let $h=t/n$. The third-order composition from Proposition~\ref{prop333} is
\[
V_3(h)
=
T_2\!\left(\frac{\overline{\alpha}h}{2}\right)
T_1(\overline{\alpha}h)
T_2\!\left(\frac{h}{2}\right)
T_1(\alpha h)
T_2\!\left(\frac{\alpha h}{2}\right).
\]
At the state level, the frozen maps are applied in the same chronological order:
\[
\Phi_2\!\left(\frac{\overline{\alpha}h}{2},\cdot\right),
~
\Phi_1(\overline{\alpha}h,\cdot),
~
\Phi_2\!\left(\frac{h}{2},\cdot\right),
~
\Phi_1(\alpha h,\cdot),
~
\Phi_2\!\left(\frac{\alpha h}{2},\cdot\right).
\]

Thus, starting from $(x_k,y_k)$, one step of the method is
\begin{align*}
\widetilde{y}
&=
\gamma_{x_k}\!\left(\frac{\overline{\alpha}h}{2},y_k\right),\\
\widetilde{x}
&=
\sigma_{\widetilde{y}}(\overline{\alpha}h,x_k),\\
\widehat{y}
&=
\gamma_{\widetilde{x}}\!\left(\frac{h}{2},\widetilde{y}\right),\\
\widehat{x}
&=
\sigma_{\widehat{y}}(\alpha h,\widetilde{x}),\\
y_{k+1}
&=
\gamma_{\widehat{x}}\!\left(\frac{\alpha h}{2},\widehat{y}\right),\\
x_{k+1}
&=
\widehat{x}.
\end{align*}
Equivalently,
\[
(x_{k+1},y_{k+1})
=
\Psi_{3,h}(x_k,y_k),
\]
where $\Psi_{3,h}$ denotes the state map defined by the five frozen substeps above. After $n$ steps, $(x_n,y_n)$ is the third-order splitting approximation at time $t$.

This formula assumes that the frozen maps $\sigma_y(s,x)$ and $\gamma_x(s,y)$ are available for the complex time parameters appearing above. In Section~\ref{sec:experiments}, this requirement is satisfied because the corresponding frozen maps have closed-form analytic extensions in the time variable.

\subsection{Three-dimensional third-order update}

We next record the three-dimensional composition used for the column labelled third order in the Lorenz experiment. This formula is obtained by applying the recursive coordinate substitution from Subsection~\ref{expHI} to the two-operator third-order composition in Proposition~\ref{prop333}. It contains $13$ frozen-flow factors, in agreement with Table~\ref{tabNterms}.

Let $\alpha=\frac{1}{2}+\frac{i\sqrt{3}}{6}$, and let $h=t/n$. The chronological order of frozen state maps is
\[
\begin{array}{r@{\;}l r@{\;}l r@{\;}l r@{\;}l r@{\;}l}
1.& \Phi_2\!\left(\frac{\overline{\alpha}h}{2}\right),&
2.& \Phi_3\!\left(\frac{\overline{\alpha}^{\,2}h}{2}\right),&
3.& \Phi_1(\overline{\alpha}^{\,2}h),&
4.& \Phi_3\!\left(\frac{\overline{\alpha}h}{2}\right),&
5.& \Phi_1(\alpha\overline{\alpha}h),\\[0.4em]
6.& \Phi_3\!\left(\frac{\alpha\overline{\alpha}h}{2}\right),&
7.& \Phi_2\!\left(\frac{h}{2}\right),&
8.& \Phi_3\!\left(\frac{\alpha\overline{\alpha}h}{2}\right),&
9.& \Phi_1(\alpha\overline{\alpha}h),&
10.& \Phi_3\!\left(\frac{\alpha h}{2}\right),\\[0.4em]
11.& \Phi_1(\alpha^2h),&
12.& \Phi_3\!\left(\frac{\alpha^2h}{2}\right),&
13.& \Phi_2\!\left(\frac{\alpha h}{2}\right).&&&&
\end{array}
\]
Here, $\Phi_1$, $\Phi_2$, and $\Phi_3$ are the frozen maps associated with the $x$, $y$, and $z$ coordinates, respectively.

Starting from $x^{(0)}=x_k, ~ y^{(0)}=y_k, ~ z^{(0)}=z_k$, one step of the method is computed as follows:
\[
\begin{aligned}
y^{(1)}
&=
\gamma_{x^{(0)},z^{(0)}}\!\left(\frac{\overline{\alpha}h}{2},y^{(0)}\right),\\
z^{(1)}
&=
\tau_{x^{(0)},y^{(1)}}\!\left(\frac{\overline{\alpha}^{\,2}h}{2},z^{(0)}\right),\\
x^{(1)}
&=
\sigma_{y^{(1)},z^{(1)}}\!\left(\overline{\alpha}^{\,2}h,x^{(0)}\right),\\
z^{(2)}
&=
\tau_{x^{(1)},y^{(1)}}\!\left(\frac{\overline{\alpha}h}{2},z^{(1)}\right),\\
x^{(2)}
&=
\sigma_{y^{(1)},z^{(2)}}\!\left(\alpha\overline{\alpha}h,x^{(1)}\right),\\
z^{(3)}
&=
\tau_{x^{(2)},y^{(1)}}\!\left(\frac{\alpha\overline{\alpha}h}{2},z^{(2)}\right),\\
y^{(2)}
&=
\gamma_{x^{(2)},z^{(3)}}\!\left(\frac{h}{2},y^{(1)}\right),\\
z^{(4)}
&=
\tau_{x^{(2)},y^{(2)}}\!\left(\frac{\alpha\overline{\alpha}h}{2},z^{(3)}\right),\\
x^{(3)}
&=
\sigma_{y^{(2)},z^{(4)}}\!\left(\alpha\overline{\alpha}h,x^{(2)}\right),\\
z^{(5)}
&=
\tau_{x^{(3)},y^{(2)}}\!\left(\frac{\alpha h}{2},z^{(4)}\right),\\
x^{(4)}
&=
\sigma_{y^{(2)},z^{(5)}}\!\left(\alpha^2h,x^{(3)}\right),\\
z^{(6)}
&=
\tau_{x^{(4)},y^{(2)}}\!\left(\frac{\alpha^2h}{2},z^{(5)}\right),\\
y^{(3)}
&=
\gamma_{x^{(4)},z^{(6)}}\!\left(\frac{\alpha h}{2},y^{(2)}\right).
\end{aligned}
\]
The updated state is then
\[
x_{k+1}=x^{(4)},
~
y_{k+1}=y^{(3)},
~
z_{k+1}=z^{(6)}.
\]
After $n$ steps, $(x_n,y_n,z_n)$ is the corresponding splitting approximation at time $t$.

This formula assumes that the frozen maps $\sigma_{y,z}$, $\gamma_{x,z}$, and $\tau_{x,y}$ can be evaluated for the complex-time parameters appearing above. In the Lorenz experiment, this is possible because the frozen maps are affine scalar flows with closed-form analytic dependence on the time parameter.